\documentclass[11pt]{article}

\usepackage[a4paper,left=2cm,right=2cm,top=1.5cm,bottom=2cm]{geometry}
\usepackage[utf8]{inputenc}
\usepackage[T1]{fontenc}
\usepackage[english]{babel}
\usepackage{amsthm}
\usepackage{amsmath}
\usepackage{amssymb}
\usepackage{mathrsfs}
\usepackage{setspace}
\usepackage{float}
\usepackage{multirow}
\usepackage{stmaryrd}
\usepackage{color}
\usepackage{wrapfig}
\usepackage{pifont}
\usepackage{array}
\usepackage[colorlinks=true,linkcolor=ocre,citecolor=ocre]{hyperref}
\usepackage{dsfont}
\usepackage{upgreek}
\usepackage{nicefrac}
\usepackage{tikz}{}
\usepackage{gensymb}
\usepackage{FiraSans}

\usepackage{graphicx}
\usepackage{caption}
\renewenvironment{proof}{\paragraph{\textsf{Proof :}}}{\begin{flushright}$\square$\end{flushright}}

\newcommand{\IE}{\mathbb{E}}
\newcommand{\IN}{\mathbb{N}}
\newcommand{\IZ}{\mathbb{Z}}

\newcommand{\IR}{\mathbb{R}}

\newcommand{\IT}{\mathbb{T}}

\newcommand{\IP}{\mathbb{P}}
\newcommand{\dst}{\displaystyle}
\newcommand{\drm}{\mathrm d}

\newcommand{\CD}{\mathcal D}

\newcommand{\CC}{\mathcal C}
\newcommand{\CH}{\mathcal H}
\newcommand{\CA}{\mathcal A}
\newcommand{\CX}{\mathcal X}

\newcommand{\CB}{\mathcal B}
\newcommand{\CG}{\mathcal G}

\newcommand{\SA}{\mathscr{A}}

\newcommand{\IDC}{\mathds{1}}

\renewcommand{\P}{\mathsf{P}}
\newcommand{\PT}{\widetilde{\P}}
\newcommand{\PI}{\mathsf{\Pi}}

\newcommand{\DD}{\mathsf{D}}
\newcommand{\DC}{\mathsf{C}}

\newcommand{\DS}{\mathsf{S}}

\newcommand{\DO}{\mathsf{A}}

\newcommand{\x}{x} 
\newcommand{\di}{v} 

\definecolor{ocre}{RGB}{64,123,121}
\definecolor{S}{rgb}{0.0,0.5,0.0}

\newcounter{item}
\numberwithin{item}{section}

\newtheorem{theorem}[item]{\sffamily Theorem}

\newtheorem{proposition}[item]{\sffamily Proposition}
\newtheorem{lemma}[item]{\sffamily Lemma}
\newtheorem{corollary}[item]{\sffamily Corollary}

\newtheorem*{theorem*}{\sffamily Theorem}
\newtheorem*{definition*}{\sffamily Definition}
\newtheorem*{proposition*}{\sffamily Proposition}
\newtheorem*{lemma*}{\sffamily Lemma}
\newtheorem*{corollary*}{\sffamily Corollary}

\usepackage[explicit]{titlesec}
\titleformat{\section}{\centering\Large\bfseries}{\thesection \ --}{0.7em}{\Large\bfseries #1}
\titleformat{\subsection}{\centering\large\bfseries}{\thesubsection \ --}{0.4em}{\large\bfseries #1}
\titleformat{\subsubsection}{\centering\bfseries}{\thesubsubsection \ --}{0.4em}{\bfseries #1}

\let\emph\relax
\DeclareTextFontCommand{\emph}{\bfseries\em}

\providecommand{\MSC}[1]
{
	{\footnotesize	
	\textbf{MSC $\mathbf{2020}$ --} #1}
}

\providecommand{\keywords}[1]
{
	{\footnotesize	
	\textbf{Keywords --} #1}
}

\title{\bfseries Weyl law for the Anderson Hamiltonian on a two-dimensional manifold}
\author{Antoine MOUZARD}
\date{}

\begin{document}

\maketitle

\abstract{We define the Anderson Hamiltonian $H$ on a two-dimensional manifold using high order paracontrolled calculus. It is a self-adjoint operator with pure point spectrum. We get lower and upper bounds on its eigenvalues which imply an almost sure Weyl-type law for $H$.}

\bigskip

\MSC{35J10; 60H25; 58J05}


\tableofcontents

\vspace{0.5cm}

\section*{Introduction}

The study of singular stochastic Partial Differential Equations (PDEs) has rapidly grown over the last decade. Following the theory of Lyons' rough paths and Gubinelli's controlled paths developed for singular stochastic ordinary differential equations, new tools have appeared to describe solutions of such PDEs that share the same philosophy. One first constructs a random space of functions/distributions from the noise through a renormalisation step; this is purely probabilistic. One then solves the PDE with classical methods on this random space; this is purely deterministic. The litterature is also growing and two different approaches have emerged. The first is based on a local description of distributions which satisfies a precise algebraic structure in order to reassemble into global objects; this is the theory of regularity structures as devised by Hairer in \cite{Hai14}. The second approach works directly with global objects and uses tools from harmonic analysis to study products; this is the paracontrolled calculus designed by Gubinelli, Imkeller and Perkowski in \cite{GIP}. In both cases, the equation dictates via a fixed point a space of solutions built from the rough source term of the PDE. There exists interesting relations between the local and the global points of views, see for example the works \cite{BH1,BH2,MP}. As far as the renormalisation step is concerned, one has to give a meaning to a number of ill-defined functionnals of the noise; this is how singular products are dealed with. If the list of such terms is given by the equation, their construction can be performed independantly of the resolution of the PDE.

\medskip

To work on manifolds, one has to adapt these methods. For the local approach, Dahlqvist, Diehl and Driver have adapted regularity structures using local charts to study the parabolic Anderson model on Riemann surfaces, see \cite{DDD}. For the global approach, Bailleul, Bernicot and Frey in \cite{BB1,BB2} used harmonic analysis tools built from the heat semigroup instead of Fourier analysis and developed paracontrolled calculus on manifolds. As in the initial work \cite{GIP} of Gubinelli, Imkeller and Perkwoski, this was only a first order calculus and it constrained the roughness one could deal with. Bailleul and Bernicot then generalised it to a high order paracontrolled calculus in \cite{BB3} and extended the range of regularity one can consider, as far as the analytical step of the problem is concerned, again working on manifolds.

\medskip

The Anderson Hamiltonian is given by
$$
H:=\Delta+\xi
$$
where $\xi$ is a space white noise. It is for example involved in the study of evolution equations such as the heat equation with random multiplicative noise
$$
\partial_tu=\Delta u+u\xi
$$
called the Parabolic Anderson model. It first appeared in \cite{Anderson} as a description of a physical phenomena involving quantum-mechanical motion with an effect of mass concentration called Anderson localization. Since then, a large class of Schrödinger operators with random potential have beed investigated, for example as a description of random dynamics in random environment. For particular models, see the book \cite{Konig} of König for a complete survey in a discrete space setting. One of the most important question on this type of operator lies in the line of the first work by Anderson, that is whether it exhibits a localisation phenomenon. This corresponds to localisation in space of eingenfunctions associated to a part of the spectrum. While this was mainly answered in the discrete setting, very accurate results were obtained by Dumaz and Labbé in \cite{DumazLabbe1} for the continuous Anderson Hamiltonian in one dimension, see also their recent works \cite{DumazLabbe2,DumazLabbe3}. The litterature about localisation of such operators is now quite vast both in mathematics and physics. Formally, the Anderson Hamiltonian with white noise potential is the continuum limit of the discrete Anderson Hamiltonian with independent and identically distrubuted random potential with finite variance. See for example the work \cite{MPdiscrete} by Martin and Perkowski.

\medskip

In dimension $1$, the noise is regular enough for the multiplication to make sense and the operator has been constructed by Fukushima and Nakao in \cite{FN} without renormalisation using Dirichlet space methods. In two dimensions using paracontrolled calculus, Allez and Chouk were the first to construct the operator on the torus, see \cite{AllezChouk}. They introduced the space of strongly paracontrolled distributions to get an operator from $L^2$ to itself with a renormalisation procedure and proved self-adjointness with pure point spectrum. They gave bounds on its eigenvalues and a tail estimate for the largest one. They also studied the large volume limit and gave a bound on the rate of divergence. Then Labbé constructed the operator in dimension $\le3$ in \cite{Labbe} with different boundary coundition using regularity structures. It relies on a reconstruction theorem in Besov spaces from his work \cite{HairerLabbe} with Hairer. He also proved self-adjointess with pure point spectrum and gave tail estimate for all the eigenvalues as well as bounds for the large volume limit. Chouk and van Zuijlen also studied the large volume limit in two dimensions, see \cite{CVZ}. Finally Gubinelli, Ugurcan and Zachhuber constructed in \cite{GUZ} the operator in dimension $2$ and $3$ on the torus using a different approach. With a finer description of the paracontrolled structure, they showed density of the domain in $L^2$ before studying the operator. They also proved self-adjointness with pure point spectrum considering the bilinear form associated to $H$ and considered evolution PDEs associated to the Anderson Hamiltonian such as the Schrödinger equation or the wave equation. Zachhuber used this approach in \cite{Zachhuber} to prove Strichartz estimate in two dimensions, the problem for $d=3$ being the use of a Hopf-Cole type transformation to construct the domain.

\medskip

To the best of our knowledge, the present work is the first to deal with the construction of the Anderson Hamiltonian on a manifold. In particular, the paracontrolled approach with the heat semigroup deals naturaly with Sobolev spaces on a manifold while we are not aware of any adaptation of the work \cite{HairerLabbe} of Hairer and Labbé in a manifold setting. We are able to recover geometric information on the manifold from the spectral properties of the Anderson Hamiltonian as one can do from the Laplacian. For example, we recover the volume of $M$ via a Weyl law from the estimates on the spectrum. This raises many interesting associated questions. As far as PDEs are concerned, it also appears in a number of singular SPDEs of interest. For example, the Schrödinger equation has been studied on manifolds by Burq, Gérard and Tzvetkov in \cite{BGT} where they prove Strichartz inequalities. Similar questions for the stochastic version of the equation are natural to ask on manifold and has been investigated in the flat case by Zachhuber in \cite{Zachhuber}. This question was investigated in \cite{MZ} by Mouzard and Zachhuber on a two-dimensional manifold with or without boundary, both for the Schrödinger and the wave equations. In particular, the result for the wave equation strongly relies on the Weyl-type law obtained in the present work, even in the flat case of the torus.

\bigskip

{\noindent\large\textbf{Main results}}

\medskip

In this work, we construct the Anderson Hamiltonian on a two-dimensional manifold using the high order paracontrolled calculus. We adapt the space-time construction \cite{BB2,BB3} of Bailleul, Bernicot and Frey to the spatial setting and work with Sobolev spaces; in particular this work is self-contained and can serve as an introduction to the work \cite{BB3} on the high order paracontrolled calculus. The simpler spatial setting forms a gentle introduction to grasp the space-time paracontrolled calculus, the only technical difficulty being the use of Sobolev spaces in addition to the Hölder spaces. We emphasize that these tools are of interest on their own in the study of singular elliptic PDEs on manifolds and somewhat flexible to use. As application, it yields existence and uniqueness to the nonlinear Schrödinger equation with multiplicative white noise on a two-dimensional manifold using a Brezis-Gallouët type inequality. In particular, this work removes the need of the  "strongly paracontrolled distributions" introduced by Allez and Chouk and used by Gubinelli, Ugurcan and Zachhuber with a second order expansion rather than an ad-hoc modification of the first order expansion.

\medskip

Given a regularisation $\xi_\varepsilon$ of the noise, we consider the renormalised Hamiltonians
$$
H_\varepsilon:=-\Delta+\xi_\varepsilon-c_\varepsilon
$$
with $c_\varepsilon:=\IE[\PI(X_\varepsilon,\xi_\varepsilon)]$ and $X_\varepsilon:=-L^{-1}\xi_\varepsilon$. While in the flat case of the torus $c_\varepsilon$ is a renormalisation constant, it is only a renormalisation function in the general setting of a two-dimensional manifold. The renormalisation of the noise yields a enhanced noise $\Xi\in\CX^\alpha$, this is done in Section \ref{subsecrenorm}. The Anderson Hamiltonian is constructed as the limit of this operator in the resolvent norm sense and we are able to get upper and lower bounds for its eigenfunctions from the eigenfunctions of the Laplacian. This is the content of Theorems \ref{SpectralResult}, \ref{EVestimates} and \ref{ResolventConvergence}.

\medskip

\begin{theorem*}
There exists a self-adjoint operator $H$ on $L^2(M)$ with pure point spectrum $\big(\lambda_n(\Xi)\big)_{n\ge1}$ such that 
$$
H=\lim_{\varepsilon\to0}-\Delta+\xi_\varepsilon-c_\varepsilon
$$
in the resolvent sense. Moreover, there exists constants $m_\delta^-(\Xi),m_\delta^+(\Xi)$ such that
$$
\lambda_n-m_\delta^-(\Xi)\le\lambda_n(\Xi)\le(1+\delta)\lambda_n+m_\delta^+(\Xi)
$$
for any $\delta\in(0,1)$ and $n\ge1$.
\end{theorem*}

\medskip

This implies the following two corollaries, that is the upper bouds on the repartition functions of the eigenfunctions and the almost sure Weyl-type law. This corresponds to corollaries \ref{TailEstimates} and \ref{weyl}.

\medskip

\begin{corollary*}
For any $n\ge1$ and $\lambda\in\IR$, we have
$$
1-me^{-h(\lambda-2\lambda_n)^{\frac{1}{12}}}\le\IP(\lambda_n(\Xi)\le\lambda)\le me^{-h(\lambda_n-\lambda)^{\frac{1}{5}}}
$$
where $m=\IE\left[e^{h\|\Xi\|_{\CX^\alpha}}\right]$ and $h>0$ a constant.
\end{corollary*}

\medskip

\begin{corollary*}
We have almost surely
$$
\lim_{\lambda\to\infty}\lambda^{-1}|\{\lambda_n(\Xi)\le\lambda\}=\frac{\textup{Vol}(M)}{4\pi}.
$$
\end{corollary*}

\bigskip

{\noindent\large\textbf{Organisation of the paper}}

\medskip

In the first section, we introduce the approximation theory based on the heat semigroup and use it to build the paracontrolled calculus. The second section is devoted to the construction and study of the Anderson Hamiltonian $H$ on a manifold in two dimensions. We show self-adjointness with pure point spectrum and provide lower and upper bounds for the eigenvalues. We finally study the cubic Schrödinger equation in Section \ref{subsecSCH}. Appendix \ref{AppendixApprox} contains all the technical details of the approximation theory and Appendix \ref{AppendixPC} gives the proof of different continuity estimates for the paracontrolled calculus.

\medskip

The main ingredients for this work are the following. We adapt the work of Bailleul and Bernicot \cite{BB2} and give the tools of the spatial high order paracontrolled calculus. These are of interest in themselves to solve other elliptic problems, on manifold or not, and are very flexible. As Gubinelli, Ugurcan and Zachhuber in \cite{GUZ}, our method relies on the almost duality property between the resonant term and the paraproduct. Finally, we introduce a truncated paraproduct $\P^s$ to describe product on adapted scales with its companion $\PT^s$ that describes associated mild formulation.

\bigskip

\textbf{Acknowledgements :} The author would like to thank Massimiliano Gubinelli for the invitation to the Haussdorf Institute of Mathematics (HIM) where he enjoyed a lot of fruitful mathematical discussions with him and others. In particular, the author is grateful to Immanuel Zachhuber for interesting discussions that motivated the present work. The author would also like to thank Nikolay Tzvetkov and Ismaël Bailleul for readings of the preliminary drafts of this version.

\bigskip

\section{Heat semigroup and paracontrolled calculus}

On a manifold $M$, the heat semigroup 
$$
P=(P_t)_{t>0}:=(e^{tL})_{t>0}
$$ 
associated to a nice enough second order differential operator $L$ can be used to regularize distributions in $\CD'(M)$ while the tools of Fourier analysis of the flat case can not be used. One can then consider the Calder\'on decomposition as an analog of the Paley-Littlewood decomposition with a continuous scaling parameter and 
$$
Q_t:=-t\partial_tP_t
$$ 
acting like a localizer on ``frequency'' of order $t^{-\frac{1}{2}}$. After giving the geometric framework, we introduce the standard families of operators we shall use to define the Besov spaces on $M$. We then construct the paraproducts $\P$ and $\PT$ with a number of tools of the high order paracontrolled calculus needed to study elliptic singular PDEs.

\bigskip

\subsection{Geometric framework}

Let $(M,d,\mu)$ be a complete volume doubling measured Riemannian manifold. We assume $M$ compact so spatial weight are not needed; everything in this section should work in the unbounded setting of \cite{BB2}. All the kernels we consider are with respect to this measure $\mu$. Let $(V_i)_{1\le i\le\di}$ be a family of smooth vector fields identified with first order differential operators on $M$. Consider the associated second order operator $L$ given by
$$
L=-\sum_{i=1}^\di V_i^2.
$$
We assume that $L$ is \emph{elliptic}. In particular, it implies that the vector fields $(V_i)_{1\le i\le\di}$ span smoothly at every point of $M$ the tangent space and the existence of smooth functions $(\gamma_i)_{1\le i\le\di}$ such that for any $f\in C^1(M,\IR)$ and $x\in M$, we have
$$
\nabla f(x)=\sum_{i=1}^\di\gamma_i(x)V_i(f)(x)V_i(x).
$$
It also implies that $L$ is sectorial in $L^2$ with kernel the constant functions, it has a bounded $H^\infty$-calculus on $L^2$ and $-L$ generates a holomorphic semigroup $(e^{-tL})_{t>0}$ on $L^2$, see \cite{DM}. Given any collection $I=(i_1,\ldots,i_n)\in\{1,\ldots,\di\}^n$, we denote by $V_I:=V_{i_n}\ldots V_{i_1}$ the differential operator of order $|I|:=n$. Under the smoothness and ellipticity conditions, the semigroup has regularity estimate at any order, that is $(t^{\frac{|I|}{2}}V_I)e^{-tL}$ and $e^{-tL}(t^{\frac{|I|}{2}}V_I)$ have kernels $K_t(x,y)$ for any $t>0$ and $x,y\in M$ that satify the Gaussian estimates
$$
\big|K_t(x,y)\big|\lesssim\mu\big(B(x,\sqrt{t})\big)^{-1}e^{-c\frac{d(x,y)^2}{t}}
$$
and for $x'\in M$
$$
\big|K_t(x,y)-K_t(x',y)\big|\lesssim\frac{d(x,x')}{\sqrt{t}}\mu\big(B(x,\sqrt{t})\big)^{-1}e^{-c\frac{d(x,y)^2}{t}}
$$
for $d(x,x')\le\sqrt{t}$ and a constant $c>0$. The range of application contains the case of a bounded domain with its Laplacian associated with periodic or Dirichlet boundary conditions if the boundary is sufficiently regular, see again \cite{DM}. Note that the Laplacian can indeed be written in the Hörmander form, see Strook's book \cite{Stroock} for example. The operator $L:\CH^2\subset L^2\to L^2$ is not invertible since its kernel contains constant function however it is invertible up to a smooth error term. Indeed, setting
$$
L^{-1}:=\int_0^1e^{-tL}\drm t,
$$
we have $L\circ L^{-1}=\textup{Id}$ up to the regularizing operator $e^{-L}$.

\bigskip

\subsection{Approximation theory}

All computations below make sense for a choice of large enough integers $b$ and $\ell$ that are fixed in any application, we also assume $b$ even. Given $x,y\in M$ and $t\in(0,1]$, we define the Gaussian kernel
$$
\CG_t(x,y):=\frac{1}{\mu\left(B(x,\sqrt{t})\right)}\left(1+c\ \frac{d(x,y)^2}{t}\right)^{-\ell}
$$
with $c>0$ a constant. We do not emphasize the dependance on the postive constant $c$ and abuse notation by writing the same letter $\CG_t$ for two functions corresponding to two different values of the constant. We have for any $s,t\in(0,1]$
$$
\int_M\CG_t(x,y)\CG_s(y,z)\drm y\lesssim\CG_{t+s}(x,z).
$$
A choice of constant $\ell$ large enough ensure that
$$
\sup_{t\in(0,1]}\sup_{x\in M}\int_M\CG_t(x,y)\drm y<\infty.
$$
This implies that any linear operator with a kernel pointwisely bounded by $\CG_t$ is bounded in $L^p(M)$ for every $p\in[1,\infty]$. The family $(\CG_t)_{t\in(0,1]}$ is our reference kernel for Gaussian operator; this is the letter `$\mathsf{G}$' in the following definition.

\medskip

\begin{definition*}
We define $\mathsf{G}$ as the set of families $(P_t)_{t\in(0,1]}$ of linear operator on $M$ with kernels pointwisely bounded by
$$
\left|K_{P_t}(x,y)\right|\lesssim\CG_t(x,y)
$$
given any $x,y\in M$.
\end{definition*}

\medskip

We consider two such families of operators $(Q_t^{(b)})_{t\in(0,1]}$ and $(P_t^{(b)})_{t\in(0,1]}$ defined as
$$
Q_t^{(b)}:=\frac{(tL)^be^{-tL}}{(b-1)!}\quad\text{and}\quad-t\partial_tP_t^{(b)}=Q_t^{(b)}
$$
with $P_0^{(b)}=\text{Id}$. In particular, there exist a polynomial $p_b$ of degree $(b-1)$ such that $P_t^{(b)}=p_b\left(tL\right)e^{-tL}$ and $p_b(0)=1$. The family $(P_t)_{t\in(0,1]}$ regularizes distributions while the family $(Q_t)_{t\in(0,1]}$ is a kind of localizer on `frequency' of order $t^{-\frac{1}{2}}$ as one can see with the parabolic scaling of the Gaussian kernel. In the flat framework of the torus, this can be explicitly written using Fourier theory. These tools also enjoy cancellation properties as Fourier projectors however it is not as precise since the operators involved here are not locally supported. For example, the following simple computation show that the composition
$$
Q_t^{(b)}\circ Q_s^{(b)}\simeq\left(\frac{ts}{(t+s)^2}\right)^bQ_{t+s}^{(2b)}
$$
is small for $s\ll t$ or $t\ll s$ but not equal to $0$. The importance of the parameter $b$ appears here as a `degree' of cancellation. One can also see that in the fact that for any polynomial function $p$ of degree less than $2b$, we have $P_t^{(b)}p=p$ and $Q_t^{(b)}p=0$ for any $t\in(0,1]$. We now define the standard family of Gaussian operators with cancellation that we shall use in this work.

\medskip

\begin{definition*}
Let $a\in\llbracket0,2b\rrbracket$. We define the standard collection of operators with cancellation of order $a$ as the set $\mathsf{StGC}^a$ of families
$$
\left((t^{\frac{|I|}{2}}V_I)(tL)^{\frac{j}{2}}P_t^{(c)}\right)_{t\in(0,1]}
$$
with $I,j$ such that $a=|I|+j$ and $c\in\llbracket1,b\rrbracket$. These operators are uniformly bounded in $L^p(M)$ for every $p\in[1,\infty]$ as functions of the parameter $t\in(0,1]$. In particular, a standard family of operator $Q\in\mathsf{StGC}^a$ can be seen as a bounded map $t\mapsto Q_t$ from $(0,1]$ to the space of bounded linear operator on $L^p(M)$. We also set
$$
\mathsf{StGC}^{[0,2b]}:=\bigcup_{0\le a\le 2b}\mathsf{StGC}^a.
$$
\end{definition*}

\medskip

Since the first order differential operators $V_i$ do not a priori commute with each other, they do not commute with $L$ and we introduce the notation
$$
\big(V_I\phi(L)\big)^\bullet:=\phi(L)V_I
$$
for any function $\phi$ such that $\phi(L)$ is defined in order to state the following cancellation property. This is not related to any notion of duality in general. In particular, $L$ is not supposed self-adjoint here.

\medskip

\begin{proposition}
Given $a,a'\in\llbracket0,2b\rrbracket$, let $Q^1\in\mathsf{StGC}^a$ and $Q^2\in\mathsf{StGC}^{a'}$. Then for any $s,t\in(0,1]$, the composition $Q_s^1\circ Q_t^{2\bullet}$ has a kernel pointwisely bounded by
\begin{align*}
\left|K_{Q_s^1\circ Q_t^{2\bullet}}(x,y)\right|&\lesssim\left(\left(\frac{s}{t}\right)^{\frac{a}{2}}\IDC_{s<t}+\left(\frac{t}{s}\right)^{\frac{a'}{2}}\IDC_{s\ge t}\right)\CG_{t+s}(x,y)\\
&\lesssim\left(\frac{ts}{(t+s)^2}\right)^{\frac{a}{2}}\CG_{t+s}(x,y)
\end{align*}
with $a=\min(a,a')$.
\end{proposition}

\medskip

\begin{proof}
Let $t\in(0,1]$. We have
$$
Q_t^1=t^{\frac{a}{2}}V_{I}L^{\frac{j}{2}}P_t^{(c)}\quad\text{and}\quad Q_t^2=t^{\frac{a'}{2}}V_{I'}L^{\frac{j'}{2}}P_t^{(c')}
$$
with $c,c'\in\llbracket1,b\rrbracket$, $a=|I|+j$ and $a'=|I'|+j'$. For any $t,s\in(0,1]$, the composition is given by
\begin{align*}
Q_s^1\circ Q_t^{2\bullet}&=s^{\frac{a}{2}}t^{\frac{a'}{2}}V_{I}L^{\frac{j+j'}{2}}P_s^{(c)}P_t^{(c')}V_{I'}\\
&=\frac{s^{\frac{a}{2}}t^{\frac{a'}{2}}}{(t+s)^{\frac{a+a'}{2}}}(t+s)^{\frac{a+a'}{2}}V_IL^{\frac{j+j'}{2}}P_s^{(c)}P_t^{(c')}V_{I'}
\end{align*}
and this yields
\begin{align*}
K_{Q_s^1\circ Q_t^{2\bullet}}(x,y)&\lesssim\frac{s^{\frac{a}{2}}t^{\frac{a'}{2}}}{(t+s)^{\frac{a+a'}{2}}}\ \CG_{t+s}(x,y)\\
&\lesssim\Big\{\left(\frac{s}{t}\right)^{\frac{a}{2}}\IDC_{s<t}+\left(\frac{t}{s}\right)^{\frac{a'}{2}}\IDC_{s\ge t}\Big\}\CG_{t+s}(x,y).
\end{align*}
The last estimate follows from a direct computation.
\end{proof}

\medskip

Operators with cancellation but not in this standard form also appear in the description of solutions to PDEs. This is the role of the set $\mathsf{GC}^a$ of the following definition.

\medskip

\begin{definition*}
Let $a\in\llbracket0,2b\rrbracket$. We define the subset $\mathsf{GC}^a\subset\mathsf{G}$ as families $(Q_t)_{t\in(0,1]}$ of operators with the following cancellation property. For any $s,t\in(0,1]$ and standard family $S\in\mathsf{StGC}^{a'}$ with $a'\in\llbracket a,2b\rrbracket$, the operator $Q_s\circ S_t^\bullet$ has a kernel pointwisely bounded by
$$
\left|K_{Q_s\circ S_t^\bullet}(x,y)\right|\lesssim\left(\frac{ts}{(t+s)^2}\right)^{\frac{a}{2}}\CG_{t+s}(x,y).
$$
\end{definition*}

\medskip

The set $\mathsf{StGC}$ can be used to define Besov spaces on a manifold. For any $f\in L^p(M)$ with $p\in[1,\infty[$ or $f\in C(M)$, we have the following reproducing Calderón formula
$$
f=\lim_{t\to0} P_t^{(b)}f=\int_0^1 Q_t^{(b)}f\frac{\drm t}{t}+P_1^{(b)}f.
$$
We interpret it as an analog to the Paley-Littlewood decomposition of $f$ on a manifold but with a continuous parameter. Indeed, the measure $\frac{\drm t}{t}$ gives unit mass to the dyadic intervals $[2^{-(i+1)},2^{-i}]$ with the operator $Q_t^{(b)}$ as a kind of multiplier roughly localized at frequencies of size $t^{-\frac{1}{2}}$. This motivates the following definition.

\medskip

\begin{definition*}
Given any $p,q\in[1,\infty]$ and $\alpha\in(-2b,2b)$, we define the Besov space $\CB_{p,q}^\alpha(M)$ as the set of distribution $f\in\CD'(M)$ such that
$$
\|f\|_{\CB_{p,q}^\alpha}:=\left\|e^{-L}f\right\|_{L^p(M)}+\sup_{\underset{|\alpha|<k\le 2b}{Q\in\mathsf{StGC}^k}}\left\|t^{-\frac{\alpha}{2}}\|Q_tf\|_{L_x^p}\right\|_{L^q(t^{-1}\drm t)}<\infty.
$$
\end{definition*}

\begin{remark}
As far as regularity is concerned, a limitation appears with this definition of $\CB_{p,q}^\alpha$ since we can only work with regularity exponent $\alpha\in(-2b,2b)$. This is only technical and $b$ can be taken as large as needed.
\end{remark}

\medskip

The Hölder spaces $\CC^\alpha:=\CB_{\infty,\infty}^\alpha$ and Sobolev spaces $\CH^\alpha:=\CB_{2,2}^\alpha$ are of particular interest with
$$
\|f\|_{\CC^\alpha}:=\|e^{-L}f\|_{L^\infty}+\sup_{\underset{|\alpha|<k\le 2b}{Q\in\mathsf{StGC}^k}}\sup_{t\in(0,1]}t^{-\frac{\alpha}{2}}\|Q_tf\|_{L_x^\infty}
$$
and
$$
\|f\|_{\CH^\alpha}:=\|e^{-L}f\|_{L^2}+\sup_{\underset{|\alpha|<k\le 2b}{Q\in\mathsf{StGC}^k}}\left(\int_0^1 t^{-\alpha}\|Q_tf\|_{L_x^2}^2\frac{\drm t}{t}\right)^{\frac{1}{2}}.
$$
This is indeed a generalisation of the classical Hölder spaces as stated in the following proposition. We shall denote $C^\alpha$ the classical spaces of Hölder functions with the norm
$$
\|f\|_{C^\alpha}:=\|f\|_{L^\infty}+\sup_{x\neq y}\frac{|f(x)-f(y)|}{d(x,y)^\alpha}
$$
for $0<\alpha<1$. Note that for any integer regularity exponent, $\CC^\alpha\neq C^\alpha$ since $\CC^1$ is the space of Lipschitz functions. The proof of the following proposition is left to the reader, it works exactly as Proposition $5$ in \cite{BB2}.

\medskip

\begin{proposition*}
For any $\alpha\in(0,1)$, we have $\CC^\alpha=C^\alpha$ and the norms $\|\cdot\|_{\CC^\alpha}$ and $\|\cdot\|_{C^\alpha}$ are equivalent. 
\end{proposition*}

\medskip

We have an analog result for Sobolev spaces however one has to be careful in the case of a manifold with boundary. The semigroup is obtained with Dirichlet conditions hence the associated Sobolev spaces are the analog of the classical $H_0^\alpha$ spaces. We keep the notation $\CH^\alpha$ but the reader should keep that in mind.

\smallskip

Given a distribution $f\in\CC^\alpha$ and $Q\in\mathsf{StGC}^k$, we have by definition a bound for $\|Q_tf\|_{\infty}$ only for $|\alpha|<k$. If $f$ is a distribution and not a function, the quantity diverges and we still have the estimate for all $k$; this will be important to keep an accurate track of the regularity. The same holds for negative Sobolev spaces.

\medskip

\begin{proposition}\label{DistributionEstimate}
Let $-2b<\alpha<0$ and $P\in\mathsf{StGC}^k$ with $k\in\llbracket0,b\rrbracket$. For $f\in\CC^\alpha$, we have
$$
\sup_{t\in(0,1]}t^{-\frac{\alpha}{2}}\|P_tf\|_{L^\infty}\lesssim\frac{1}{k-\alpha}\|f\|_{\CC^\alpha}.
$$
For $f\in\CH^\alpha$, we have
$$
\|t^{-\frac{\alpha}{2}}\|P_tf\|_{L_x^2}\|_{L^2(t^{-1}\drm t)}\lesssim\frac{1}{k-\alpha}\|f\|_{\CH^\alpha}.
$$
\end{proposition}

\medskip

\begin{proof}
Since $P\in\mathsf{StGC}^k$ with $k\in\llbracket0,2b\rrbracket$, there exist $I=(i_1,\ldots,i_n),j\in\IN$ and $c\in\llbracket1,b\rrbracket$ such that $k=|I|+j$ and
$$
P_t=(t^{\frac{|I|}{2}}V_I)(tL)^{\frac{j}{2}}P_t^{(c)}.
$$
If $|\alpha|<k$, the result holds by definition of $\CC^\alpha$. If $|\alpha|\ge k$, we have
\begin{align*}
P_tf&=(t^{\frac{|I|}{2}}V_I)(tL)^{\frac{j}{2}}\left(\int_t^1Q_s^{(c)}f\frac{\drm s}{s}+P_1^{(c)}f\right)\\
&=\int_t^1\left(\frac{t}{s}\right)^{\frac{k}{2}}(s^{\frac{|I|}{2}}V_I)(sL)^{\frac{j+c}{2}}P_s^{(1)}f\frac{\drm s}{s}+R_tf\\
&=\int_t^1\left(\frac{t}{s}\right)^{\frac{k}{2}}Q_sf\frac{\drm s}{s}+R_tf
\end{align*}
with $Q_s:=(s^{\frac{|I|}{2}}V_I)(sL)^{\frac{j+c}{2}}P_s^{(1)}\in\mathsf{StGC}^{k+c}$ and $R_t:=(t^{\frac{|I|}{2}}V_I)(tL)^{\frac{j}{2}}P_1^{(c)}$. The term $R_tf$ is bounded because of the smoothing operator $P_1^{(c)}$. Since $c\ge1$, $Q$ belongs at least to $\mathsf{StGC}^{k+1}$ hence if $|\alpha|<k+1$ we have
\begin{align*}
t^{-\frac{\alpha}{2}}\|P_tf\|_{L^\infty}&\le t^{-\frac{\alpha}{2}}\int_t^1\left(\frac{t}{s}\right)^{\frac{k}{2}}\|Q_sf\|_{L^\infty}\frac{\drm s}{s}\\
&\le\|f\|_{\CC^\alpha}\int_t^1\left(\frac{t}{s}\right)^{\frac{k-\alpha}{2}}\frac{\drm s}{s}\\
&\le\|f\|_{\CC^\alpha}\frac{2}{k-\alpha}
\end{align*}
and this yields the result using that $\alpha<0\le k$ hence $k-\alpha>0$. If $|\alpha|\ge k+1$, using the same integral representation for $Q$ and an induction completes the proof of the $L^\infty$-estimate. For the $L^2$-estimate, we interpolate between $L^1$ and $L^\infty$ as in Appendix \ref{AppendixApprox} to get
\begin{align*}
\|t^{-\frac{\alpha}{2}}\|P_tf\|_{L^2}\|_{L^2(t^{-1}\drm t)}&\le\left\|t^{-\frac{\alpha}{2}}\int_t^1\left(\frac{t}{s}\right)^{\frac{k}{2}}\|Q_sf\|_{L^2}\frac{\drm s}{s}\right\|_{L^2(t^{-1}\drm t)}\\
&\le\frac{2}{k-\alpha}\|f\|_{\CH^\alpha}.
\end{align*}
\end{proof}

\medskip

One can see that the bound diverges as $\alpha$ goes to $0$ if the operator does not encode any cancellation, that is $k=0$. In the case $\alpha=0$, we have $\|P_tf\|_{L^\infty}\lesssim\|f\|_{L^\infty}$ hence the $L^\infty$-bound holds. However the $L^2$-bound is not satisfied since $\|P_tf\|_{L^2}\lesssim\|f\|_{L^2}$ only implies
$$
\int_0^1\|P_tf\|_{L^2}^2\frac{\drm t}{t}\le\|f\|_{L^2}^2\int_0^1\frac{\drm t}{t}=\infty.
$$
This will explain an important difference for paraproducts on negative Hölder and Sobolev spaces as one can see with Propositions \ref{ParaHolder} and \ref{ParaSobolev}.

\bigskip

\subsection{Intertwined paraproducts}

We use the standard family of Gaussian operators to study the product of distributions as one can do using the Paley-Littlewood decomposition in the flat case; this lead to the definition of the paraproduct $\P$ and the resonant term $\PI$ that describe products. Then we introduce the paraproduct $\PT$ intertwined with $\P$ to describe solutions of elliptic PDEs.

\medskip

\subsubsection{Paraproduct and resonant term}

One can define the product of a distributions $f\in\CD'(M)$ with a smooth function $g\in\CD(M)$. If however the distribution $f$ belongs to a Hölder space $\CC^\alpha$ with $\alpha<0$, one might hope to do better. It is indeed the case as we can see with the next theorem which is nothing more than Young's integration condition.

\medskip

\begin{theorem*}
The multiplication $(f,g)\mapsto fg$ extends in a unique bilinear operator from $\CC^\alpha\times\CC^\beta$ to $\CC^{\alpha\wedge\beta}$ if and only if $\alpha+\beta>0$.
\end{theorem*}

\medskip

We are however interested in the case $\alpha+\beta<0$ when dealing with singular stochastic PDEs, as we are interested to stochastic ODEs where Young's condition is not verified. Following \cite{GIP}, Bailleul, Bernicot and Frey in \cite{BB1,BB2,BB3} have defined two bilinear operators $\P_fg$ and $\PI(f,g)$ such that we have the formal decomposition of the product of two distributions as
$$
fg=\P_fg+\PI(f,g)+\P_gf
$$
where the paraproducts $\P_fg$ and $\P_gf$ are well-defined for any distibutions $f,g\in\CD'(M)$. Of course, this means that $\PI(f,g)$ does have a meaning for $f\in\CC^\alpha$ and $g\in\CC^\beta$ if and only if $\alpha+\beta>0$; this is the resonant term. We want this decomposition to keep an accurate track of the regularity of each terms. More precisely, $\P_fg$ and $\PI(f,g)$ should belong to $\CC^{\alpha+\beta}$ if $\alpha<0$ while $\P_gf$ to the less regular space $\CC^\alpha$ as it is the case for the torus. We construct in this work such paraproduct and resonant term for space distributions on our manifold $M$, we mainly follow \cite{BB2} in the simpler spatial setting.

\bigskip

Let $f,g\in\CD'(M)$. Formally, we have
\begin{align*}
fg&=\lim_{t\to 0}P_t^{(b)}\left(P_t^{(b)}f\cdot P_t^{(b)}g\right)\\
&=\int_0^1\left\{Q_t^{(b)}\big(P_t^{(b)}f\cdot P_t^{(b)}g\big)+P_t^{(b)}\big(Q_t^{(b)}f\cdot P_t^{(b)}g\big)+P_t^{(b)}\big(P_t^{(b)}f\cdot Q_t^{(b)}g\big)\right\}\frac{\drm t}{t}\\
&\quad\quad\quad+P_1^{(b)}\left(P_1^{(b)}f\cdot P_1^{(b)}g\right).
\end{align*}
The last term being smooth, it does not bother us. Remark that the choice of the constant ``$1$'' is arbitrary and it might be useful to change it, as one can see with the construction of the Anderson Hamiltonian. The family $P^{(b)}$ does not encode any cancellation while $Q^{(b)}$ encodes cancellation of order $2b$ so each terms in the integral have one operator with a lot of cancellations and two with none. Since we do not have nice estimates for these terms, we want to transfer some of the cancellation from $Q^{(b)}$ to the $P^{(b)}$ in each term. To do so, we use the Leibnitz rule
$$
V_i(fg)=V_i(f)g+fV_i(g).
$$
For example, we have
\begin{align*}
\int_0^1P_t^{(b)}\left((tV_i^2)Q_t^{(b-1)}f\cdot P_t^{(b)}g\right)\frac{\drm t}{t}=&\int_0^1P_t^{(b)}(\sqrt{t}V_i)\left((\sqrt{t}V_i)Q_t^{(b-1)}f\cdot P_t^{(b)}g\right)\frac{\drm t}{t}\\
&\quad-\int_0^1P_t^{(b)}\left((\sqrt{t}V_i)Q_t^{(b-1)}f\cdot(\sqrt{t}V_i)P_t^{(b)}g\right)\frac{\drm t}{t}
\end{align*}
so if we denote by $(c_1,c_2,c_3)$ the cancellation of the three operators in the integral, we have
$$
(0,2b,0)=(1,2b-1,0)+(0,2b-1,1).
$$
This shows that we will not be able to have cancellation for all three operators at the same time but at least two. This is where the notation $Q^\bullet$ comes into play and multiple uses of this trick allows to decompose the product as
$$
fg=\sum_{\mathbf{a}\in\CA_b}\sum_{\mathbf{Q}\in\mathsf{StGC}^{\bf a}}b_{\bf Q}\int_0^1Q_t^{1\bullet}\left(Q_t^2f\cdot Q_t^3g\right)\frac{\drm t}{t}
$$
where $\mathbf{Q}=(Q_1,Q_2,Q_3)$, $\mathsf{StGC}^{\bf a}=\mathsf{StGC}^{a_1}\times\mathsf{StGC}^{a_2}\times\mathsf{StGC}^{a_3}$,
$$
\CA_b=\big\{(a_1,a_2,a_3)\in\IN^3\ ;\ a_1+a_2+a_3=2b\quad\text{and}\quad a_1,a_2\text{ or }a_3=b\big\}
$$
and $b_{\bf Q}\in\IR$ is a real coefficient associated to $\mathbf{Q}$. In particular, only one of the $a_i$ in $\mathbf{a}\in\SA_b$ can be less than $\frac{b}{2}$ and this gives us three terms $\P_fg,\P_gf$ and $\PI(f,g)$ such that
$$
fg=\P_fg+\PI(f,g)+\P_gf+P_1^{(b)}\left(P_1^{(b)}f\cdot P_1^{(b)}g\right).
$$

\medskip

\begin{definition*}
Given two distributions $f,g\in\CD'(M)$, we define the paraproduct and the resonant term as
$$
\P_fg:=\sum_{\mathbf{a}\in\CA_b;a_2<\frac{b}{2}}\sum_{\mathbf{Q}\in\mathsf{StGC}^{\bf a}}b_{\bf Q}\int_0^1 Q_t^{1\bullet}\left(Q_t^2f\cdot Q_t^3g\right)\frac{\drm t}{t}.
$$
and
$$
\PI(f,g):=\sum_{\mathbf{a}\in\CA_b;a_2,a_3\ge\frac{b}{2}}\sum_{\mathbf{Q}\in\mathsf{StGC}^{\bf a}}b_{\bf Q}\int_0^1 Q_t^{1\bullet}\left(Q_t^2f\cdot Q_t^3g\right)\frac{\drm t}{t}.
$$
\end{definition*}

\medskip

In particular, $\P_fg$ is a linear combination of
$$
\int_0^1Q_t^{1\bullet}\left(P_tf\cdot Q_t^2g\right)\frac{\drm t}{t}
$$
and $\PI(f,g)$ of
$$
\int_0^1P_t^{\bullet}\left(Q_t^1f\cdot Q_t^2g\right)\frac{\drm t}{t}
$$
with $Q^1,Q^2\in\mathsf{StGC}^{\frac{b}{2}}$ and $P\in\mathsf{StGC}^{[0,b]}$. We insist that in the following $P$ will denote an operator with possibly no cancellations while $Q$ will denote an operator with cancellations of order at least $\frac{b}{2}$.

\medskip

These operators enjoy the same continuity estimates as their Fourier counterparts from which one can recover Young's condition as stated in the following Proposition.

\medskip

\begin{proposition}\label{ParaHolder}
Let $\alpha,\beta\in(-2b,2b)$ be regularity exponents.
\begin{itemize}
	\item[$\centerdot$] If $\alpha\ge0$, then $(f,g)\mapsto\P_fg$ is continuous from $\CC^\alpha\times\CC^\beta$ to $\CC^\beta$.
	\item[$\centerdot$] If $\alpha<0$, then $(f,g)\mapsto\P_fg$ is continuous from $\CC^\alpha\times\CC^\beta$ to $\CC^{\alpha+\beta}$.
	\item[$\centerdot$] If $\alpha+\beta>0$, then $(f,g)\mapsto\PI(f,g)$ is continuous from $\CC^\alpha\times\CC^\beta$ to $\CC^{\alpha+\beta}$.
\end{itemize}
\end{proposition}

\medskip

\begin{proof}
Let us first consider the case $\alpha<0$ and let $Q\in\mathsf{StGC}^r$ with $r>|\alpha+\beta|$. Recall that $\P_fg$ is a linear combination of terms of the form
$$
\int_0^1 Q_t^{1\bullet}\left(P_tf\cdot Q_t^2g\right)\frac{\drm t}{t}
$$
with $Q^1,Q^2\in\mathsf{StGC}^{\frac{b}{2}}$ and $P\in\mathsf{StGC}^{[0,b]}$. Since $\alpha<0$, \ref{DistributionEstimate} gives
\begin{align*}
\left|\int_0^1 Q_sQ_t^{1\bullet}\left(P_tf\cdot Q_t^2g\right)\frac{\drm t}{t}\right|&\lesssim\int_0^1\left(\frac{ts}{(t+s)^2}\right)^{\frac{r}{2}}\|f\|_{\CC^\alpha}\|g\|_{\CC^\beta}t^{\frac{\alpha+\beta}{2}}\frac{\drm t}{t}\\
&\lesssim s^{\frac{\alpha+\beta}{2}}\|f\|_{\CC^\alpha}\|g\|_{\CC^\beta}
\end{align*}
for any $s\in(0,1)$ hence $\P_fg\in\CC^{\alpha+\beta}$.

\smallskip

For $\alpha\ge0$, we consider $Q\in\mathsf{StGC}^r$ with $r>|\beta|$. In this case, we have $|P_tf|\le\|f\|_{\CC^\alpha}$ for all $t\in(0,1)$ so
$$
\left|\int_0^1 Q_sQ_t^{1\bullet}\left(P_tf\cdot Q_t^2g\right)\frac{\drm t}{t}\right|\lesssim s^{\frac{\beta}{2}}\|f\|_{\CC^\alpha}\|g\|_{\CC^\beta}
$$
hence $\P_fg\in\CC^\beta$.

\smallskip

For the resonant term, let $Q\in\mathsf{StGC}^r$ with $r>|\alpha+\beta|$. We have
\begin{align*}
\left|\int_0^1 Q_sP_t^\bullet\left(Q_t^1f\cdot Q_t^2g\right)\frac{\drm t}{t}\right|&\lesssim
\|f\|_{\CC^\alpha}\|g\|_{\CC^\beta}\left(\int_0^st^{\frac{\alpha+\beta}{2}}\frac{\drm t}{t}+\int_s^1\left(\frac{s}{t}\right)^{\frac{r}{2}}t^{\frac{\alpha+\beta}{2}}\frac{\drm t}{t}\right)\\
&\lesssim s^{\frac{\alpha+\beta}{2}}\|f\|_{\CC^\alpha}\|f\|_{\CC^\beta}
\end{align*}
using that $\alpha+\beta>0$ hence $\PI(f,g)\in\CC^{\alpha+\beta}$.
\end{proof}

\medskip

We also have estimates for the Sobolev spaces whose proofs are given in Proposition \ref{ParaSobolevAppendix} from Appendix \ref{AppendixPC}.

\medskip

\begin{proposition}\label{ParaSobolev}
Let $\alpha,\beta\in(-2b,2b)$ be regularity exponents.
\begin{itemize}
	\item[$\centerdot$] If $\alpha>0$, then $(f,g)\mapsto\P_fg$ is continuous from $\CC^\alpha\times\CH^\beta$ to $\CH^\beta$ and from $\CH^\alpha\times\CC^\beta$ to $\CH^\beta$.
	\item[$\centerdot$] If $\alpha<0$, then $(f,g)\mapsto\P_fg$ is continuous from $\CC^\alpha\times\CH^\beta$ to $\CH^{\alpha+\beta}$ and from $\CH^\alpha\times\CC^\beta$ to $\CH^{\alpha+\beta}$.
	\item[$\centerdot$] If $\alpha+\beta>0$, then $(f,g)\mapsto\PI(f,g)$ is continuous from $\CH^\alpha\times\CC^\beta$ to $\CH^{\alpha+\beta}$.
\end{itemize}
\end{proposition}

\medskip

In particular, this implies that $(f,g)\mapsto\P_fg$ is continuous from $L^2\times\CC^\beta$ to $\CH^{\beta-\delta}$ for all $\delta>0$. For Sobolev spaces, there is a small loss of regularity and one does not recover the space $\CH^\beta$ while this does not happen for Hölder spaces. This comes from the remark following Proposition \ref{DistributionEstimate}.

\medskip

As in the works \cite{GUZ,Zachhuber} of Gubinelli, Ugurcan and Zachhuber, one last property of $\P$ and $\PI$ in terms of Sobolev spaces is that $\P$ is almost the adjoint of $\PI$ when $L$ is self-adjoint in the sense that the difference is more regular. A careful track of the previous computation show that for all $\mathbf{a}\in\{(0,b,b),(b,0,b),(b,b,0)\}$ and $\mathbf{Q}\in\mathsf{StGC}^{\bf a}$, we have $b_{\bf Q}=0$ except for ${\bf Q}=(P_t^{(b)},Q_t^{(b/2)},Q_t^{(b/2)}),(Q_t^{(b/2)},P_t^{(b)},Q_t^{(b/2)})$ or $(Q_t^{(b/2)},Q_t^{(b/2)},P_t^{(b)})$ where $b_{\bf Q}=1$. Define the corrector for almost duality as
$$
\DO(a,b,c):=\big\langle a,\PI(b,c)\big\rangle-\big\langle\P_ab,c\big\rangle.
$$

\medskip

\begin{proposition}\label{Duality}
Assume $L$ self-adjoint. Let $\alpha,\beta,\gamma\in(-2b,2b)$ such that $\beta+\gamma<1$ and $\alpha+\beta+\gamma\ge0$. If $\alpha<1$, then $(a,b,c)\mapsto\DO(a,b,c)$ extends in a unique trilinear operator from $\CH^\alpha\times\CC^\beta\times\CH^\gamma$ to $\IR$.
\end{proposition}

\medskip

\begin{proof}
$\DO(a,b,c)$ is a linear combination of
$$
\int_0^1\Big\{\Big\langle a,P_t^{1\bullet}\big(Q_t^1b\cdot Q_t^2c\big)\Big\rangle-\Big\langle Q_t^{3\bullet}\big(P_t^2a\cdot Q_t^4b\big),c\Big\rangle\Big\}\frac{\drm t}{t}
$$
with $P^1,P^2\in\mathsf{StGC}^{[0,b]}$ and $Q^1,Q^2,Q^3,Q^4\in\mathsf{StGC}^{\frac{b}{2}}$. We first consider $P^1,P^2\in\mathsf{StGC}^0$. By construction of the paraproduct and the resonant term, we have $P^1=P^2=P^{(b)}=:P$ and $Q^1=Q^2=Q^3=Q^4=Q^{(b/2)}=:Q$ hence we consider
$$
\int_0^1\Big\{\Big\langle a,P_t\big(Q_tb\cdot Q_tc\big)\Big\rangle-\Big\langle Q_t\big(P_ta\cdot Q_tb\big),c\Big\rangle\Big\}\frac{\drm t}{t}.
$$
Since $L$ is self-adjoint, $P_t$ and $Q_t$ are too and we have
\begin{align*}
\int_0^1\Big\langle a,P_t\big(Q_tb\cdot Q_tc\big)\Big\rangle\frac{\drm t}{t}&=\int_0^1\Big\langle P_ta,Q_tb\cdot Q_tc\Big\rangle\frac{\drm t}{t}\\
&=\int_0^1\Big\langle P_ta\cdot Q_tb,Q_tc\Big\rangle\frac{\drm t}{t}\\
&=\int_0^1\Big\langle Q_t\big(P_ta\cdot Q_tb\big),c\Big\rangle\frac{\drm t}{t}
\end{align*}
hence the difference is equal to $0$. Let us now consider the terms with $P^1,P^2\in\mathsf{StGC}^{[1,b]}$ and bound each of them independently. Since $\alpha+\beta+\gamma\ge0$, we have
\begin{align*}
\left|\int_0^1\Big\langle a,P_t^{1\bullet}\big(Q_t^2b\cdot Q_t^3c\big)\Big\rangle\frac{\drm t}{t}\right|&\lesssim\|a\|_{\CH^\alpha}\left\|\int_0^1P_t^{1\bullet}\big(Q_t^2b\cdot Q_t^3c\big)\frac{\drm t}{t}\right\|_{\CH^{\beta+\gamma}}\\
&\lesssim\|a\|_{\CH^\alpha}\|b\|_{\CC^\beta}\|c\|_{\CH^\gamma}
\end{align*}
with $\beta+\gamma<1$ and using $\alpha\in(0,1)$ we have
\begin{align*}
\left|\int_0^1\Big\langle Q_t^{3\bullet}\big(P_t^2a\cdot Q_t^4b\big),c\Big\rangle\frac{\drm t}{t}\right|&\lesssim\left\|\int_0^1Q_t^{3\bullet}\big(P_t^2a\cdot Q_t^4b\big)\frac{\drm t}{t}\right\|_{\CH^{\alpha+\beta}}\|c\|_{\CH^\gamma}\\
&\lesssim\|a\|_{\CH^\alpha}\|b\|_{\CC^\beta}\|c\|_{\CH^\gamma}
\end{align*}
which completes the proof since $\alpha+\beta+\gamma\ge0$.
\end{proof}

\subsubsection{Intertwined paraproducts}

The description of solution to elliptic PDEs involving $L$ using paracontrolled calculus necessitate to study how $L$ and $\P$ interacte with each other. Following Bailleul, Bernicot and Frey in \cite{BB2}, we want to define a new paraproduct $\PT$ intertwined with the paraproduct through
$$
L\PT_fg=\P_fLg.
$$
Since $L$ is not invertible, we use $L^{-1}$ an inverse up to a smooth error term. Hence a more conceivable intertwining relation is
$$
L\PT_fg=\P_fLg-e^{-L}\left(\P_fLg\right).
$$

\medskip

\begin{definition*}
Given any distributions $f,g\in\CD'(M)$, we define $\PT_fg$ as
$$
\PT_fg:=L^{-1}\P_fLg
$$
for which we have the explicit formula
$$
\PT_fg=\sum_{\mathbf{a}\in\SA_b;a_2<\frac{b}{2}}\sum_{\mathbf{Q}\in\mathsf{StGC}^{\bf a}}b_{\bf Q}\int_0^1\widetilde Q_t^{1\bullet}\left(Q_t^2f\cdot\widetilde Q_t^3g\right)\frac{\drm t}{t}
$$
where $\widetilde Q_t^1:=Q_t^1(tL)^{-1}$ and $\widetilde Q_t^3:=Q_t^3(tL)$.
\end{definition*}

\medskip

It is immediate that $\widetilde Q^3$ belongs to $\mathsf{StGC}^{a_3+2}$. The cancellation property of $\widetilde Q^1$ is given by the following lemma. Remark that it is not in standard form anymore, this is where the $\mathsf{GC}$ class comes into play.

\medskip

\begin{lemma}\label{lemmaPT}
Let $Q\in\mathsf{StGC}^{\frac{b}{2}}$. Then $\widetilde Q_t:=Q_t(tL)^{-1}$ defines a family that belongs to $\mathsf{GC}^{\frac{b}{2}-2}$ for $b$ large enough.
\end{lemma}

\medskip

\begin{proof}
Since $Q\in\mathsf{StGC}^{\frac{b}{2}}$, there exist $I=(i_1,\ldots,i_n),j\in\IN$ and $c\in\llbracket1,b\rrbracket$ such that $\frac{b}{2}=|I|+j$ and
$$
Q_t=(t^{\frac{|I|}{2}}V_I)(tL)^{\frac{j}{2}}P_t^{(c)}.
$$ 
This immediatly follows from
\begin{align*}
Q_t(tL)^{-1}&=(t^{\frac{|I|}{2}}V_I)(tL)^{\frac{j}{2}}(tL)^{-1}P_t^{(c)}\\
&=(t^{\frac{|I|}{2}}V_I)(tL)^{\frac{j-2}{2}}P_t^{(c)}(\textup{Id}-e^{L}).
\end{align*}
\end{proof}

\medskip

This lemma immediatly yields the following proposition, that is $\PT$ has the same structure as $\P$ hence the same continuity estimates.

\medskip

\begin{proposition}
For any distribution $f,g\in\CD'(M)$, $\PT_fg$ is given as a linear combination of terms of the form
$$
\int_0^1\widetilde Q_t^{1\bullet}\left(Q_t^2f\cdot\widetilde Q_t^3g\right)\frac{\drm t}{t}
$$
where $\widetilde Q^1\in\mathsf{GC}^{\frac{b}{2}-2},Q^2\in\mathsf{StGC}^{[0,b]}$ and $\widetilde Q^3\in\mathsf{StGC}^{\frac{b}{2}+2}$. Thus for any regularity exponent $\alpha,\beta\in(-2b,2b)$, we have the following continuity results.
\begin{itemize}
	\item[$\centerdot$] If $\alpha\ge0$, then $(f,g)\mapsto\PT_fg$ is continuous from $\CC^\alpha\times\CC^\beta$ to $\CC^\beta$.
	\item[$\centerdot$] If $\alpha<0$, then $(f,g)\mapsto\PT_fg$ is continuous from $\CC^\alpha\times\CC^\beta$ to $\CC^{\alpha+\beta}$.
\end{itemize}
We also have the same associated Sobolev estimates.
\begin{itemize}
	\item[$\centerdot$] If $\alpha>0$, then $(f,g)\mapsto\PT_fg$ is continuous from $\CC^\alpha\times\CH^\beta$ to $\CH^\beta$ and from $\CH^\alpha\times\CC^\beta$ to $\CH^\beta$.
	\item[$\centerdot$] If $\alpha<0$, then $(f,g)\mapsto\PT_fg$ is continuous from $\CC^\alpha\times\CH^\beta$ to $\CH^{\alpha+\beta}$ and from $\CH^\alpha\times\CC^\beta$ to $\CH^{\alpha+\beta}$.
\end{itemize}
\end{proposition}

\bigskip

\subsection{Correctors and commutators}

The study of elliptic PDEs with singular product involves resonant term given a function $u$ paracontrolled by a noise dependent function $X\in\CC^\alpha$, that is
$$
u=\PT_{u'}X+u^\sharp
$$
with $u'\in\CC^\alpha$ and $u^\sharp\in\CC^{2\alpha}$ a smoother remainder. If $\alpha<1$, the product $u\zeta$ is singular for $\zeta\in\CC^{\alpha-2}$ however we have the formal decomposition
$$
\PI(u,\zeta)=\PI\big(\PT_{u'}X,\zeta\big)+\PI(u^\sharp,\zeta)=u'\PI(X,\zeta)+\DC(u',X,\zeta)+\PI(u^\sharp,\zeta)
$$
with the corrector $\DC$ introduced by Gubinelli, Imkeller and Perkowski in \cite{GIP} defined as
$$
\DC(a_1,a_2,b):=\PI\big(\PT_{a_1}a_2,b\big)-a_1\PI(a_2,b).
$$
If $\frac{2}{3}<\alpha<1$, then the product $\PI(u^\sharp,\zeta)$ is well-defined. Thus we are able to give a meaning to the product $u\zeta$ for $u$ paracontrolled by $X$ once we have a proper continuity estimate for $\DC$ and a meaning to the product $X\zeta$; this is the controlled rough path philosophy. This last task is only a probabilistic one and does not impact the analytical resolution of the equation, this is the renormalisation step. We state here a continuity estimate for $\DC$ while its proof is given in Proposition \ref{CorrectorProof} in Appendix \ref{AppendixPC}. 

\medskip

\begin{proposition}\label{CorrectorHolder}
Let $\alpha_1\in(0,1)$ and $\alpha_2,\beta\in\IR$. If
$$
\alpha_2+\beta<0\quad\text{and}\quad\alpha_1+\alpha_2+\beta>0,
$$
then $(a_1,a_2,b)\mapsto\DC(a_1,a_2,b)$ extends in a unique continuous operator from $\CC^{\alpha_1}\times\CC^{\alpha_2}\times\CC^\beta$ to $\CC^{\alpha_1+\alpha_2+\beta}$.
\end{proposition}

\medskip

We also have the following proposition to work with Sobolev spaces.

\medskip

\begin{proposition}\label{CorrectorSobolev}
Let $\alpha_1\in(0,1)$ and $\alpha_2,\beta\in\IR$. If
$$
\alpha_2+\beta<0\quad\text{and}\quad\alpha_1+\alpha_2+\beta>0,
$$
then $(a_1,a_2,b)\mapsto\DC(a_1,a_2,b)$ extends in a unique continuous operator from $\CH^{\alpha_1}\times\CC^{\alpha_2}\times\CC^\beta$ to $\CH^{\alpha_1+\alpha_2+\beta}$.
\end{proposition}

\medskip

\begin{remark}
Note that the first paramater $\alpha_1$ has to be smaller than $1$. This is due to the fact that for any function $f\in\CC^\alpha$ with $\alpha\ge0$, one has
$$
|f(x)-f(y)|\le\|f\|_{\CC^\alpha}d(x,y)^{\alpha\wedge1}
$$
with a factor no greater than $1$ even if $\alpha$ is. This means that we are not able to benefit from regularity greater than $1$ with only a first order Taylor expansion. To work with a function of regularity $\alpha_1\in(1,2)$, one have to consider the refined corrector defined in the flat one dimensional case by
$$
\DC^{(1)}\big(a_1,a_2,b\big)(x):=\PI\big(\PT_{a_1}a_2,b\big)(x)-a_1(x)\PI\big(a_2,b\big)(x)-a_1'(x)\PI\big(\PT_{(x-\cdot)}a_2,b\big)(x)
$$
that we interpret as a first order refined corrector for $x\in\IT$. There is an analog refined corrector on a manifold $M$, see \cite{BB3}. However, this will not be needed in this work.
\end{remark}

\medskip

We need the corrector $\DC$ to study ill-defined product, this is the condition $\alpha_2+\beta<0$. However, we also have to investigate well-defined product to get more accurate descriptions. For this purpose, we introduce the commutator
$$
\DD(a_1,a_2,b):=\PI\big(\PT_{a_1}a_2,b)-\P_{a_1}\PI(a_2,b).
$$

\medskip

\begin{proposition}\label{CommutatorHolder}
Let $\alpha_1\in(0,1)$ and $\alpha_2,\beta\ge0$. Then $(a_1,a_2,b)\mapsto\DD(a_1,a_2,b)$ extends in a unique continuous operator from $\CC^{\alpha_1}\times\CC^{\alpha_2}\times\CC^\beta$ to $\CC^{\alpha_1+\alpha_2+\beta}$ and from $\CH^{\alpha_1}\times\CC^{\alpha_2}\times\CC^\beta$ to $\CH^{\alpha_1+\alpha_2+\beta}$.
\end{proposition}

\medskip

Again, one can bypass the condition $\alpha_1\in(0,1)$ using refined commutators. Note that in their initial work \cite{GIP}, Gubinelli, Imkeller and Perkowski call $\DC$ a commutator whereas with the point of view of high order paracontrolled calculus of \cite{BB3}, the operator $\DD$ is closer to be a commutator than $\DC$. We need one final commutator that swaps paraproducts defined by
$$
\DS(a_1,a_2,b):=\P_b\PT_{a_1}a_2-\P_{a_1}\P_ba_2.
$$

\medskip

\begin{proposition}\label{SwapHolder}
Let $\alpha_1,\alpha_2\in\IR$ and $\beta<0$. Then $(a_1,a_2,b)\mapsto\DS(a_1,a_2,b)$ extends in a unique continuous operator from $\CC^{\alpha_1}\times\CC^{\alpha_2}\times\CC^\beta$ to $\CC^{\alpha_1+\alpha_2+\beta}$ and from $\CH^{\alpha_1}\times\CC^{\alpha_2}\times\CC^\beta$ to $\CH^{\alpha_1+\alpha_2+\beta}$.
\end{proposition}

\bigskip

\section{The Anderson Hamiltonian}

In this section, we define and study the Anderson Hamiltonian
$$
H:=L+\xi
$$
where $-L$ is the Laplace-Beltrami operator on a compact two-dimensional manifold $M$ without boundary or with a smooth boundary and Dirichlet conditions. To apply the construction of the first section, one needs to have an Hörmander representation for $L$. This is possible in this case with a number of vector fields possibly greater than the dimension, see for example Section $4.2.1$ from Stroock's book \cite{Stroock}. The random potential $\xi$ is a spatial white noise and belongs almost surely to $\CC^{\alpha-2}$ for any $\alpha<1$. For a generic function $u\in L^2$, the product $u\xi$ is ill-defined hence one needs to find a proper domain for the operator. A natural method would be to take the closure of the subspace of smooth functions with the domain norm $\|u\|_{L^2}+\|Hu\|_{L^2}$. However this yields a trivial domain since $Hu$ has the same regularity as the noise, because of the product $u\xi$ if $u$ is smooth, thus it does not belong to $L^2$. Following the recent study of singular SPDEs, one can construct a random domain $\CD_\Xi$ depending on an enhancement $\Xi$ of the noise obtained through a renormalisation procedure. One can use the paraproduct to decompose the product for $u\in\CH^\alpha$ as
$$
u\xi=\P_u\xi+\P_\xi u+\PI(u,\xi).
$$
In this expression, the roughest term is $\P_u\xi\in\CC^{\alpha-2}$ while $\P_\xi u+\PI(u,\xi)$ formally belongs to $\CH^{2\alpha-2}$. For a function $u$ in the domain, we want to cancel out the roughest part of the product using the Laplacian term $Lu$, hence we want
$$
Lu=\P_u\xi+v^\sharp
$$
with $v^\sharp\in\CH^{2\alpha-2}$. This suggests the paracontrolled expansion
$$
u=\PT_uX+u^\sharp
$$
with 
$$
X:=L^{-1}\xi
$$ 
and $u^\sharp\in\CH^{2\alpha}$. We insist that we want functions in the domain to encode exactly what is needed to have a cancellation between the Laplacian and the product. In particular, $H$ is not treated at all like a perturbation of the Laplacian. 

\smallskip

At this point, two natural questions arise. Is the subspace of such paracontrolled functions dense in $L^2$ and can one make sense of the singular product? 

\begin{enumerate}
	\item[$1)$] For the first question, one can introduce a parameter $s>0$, in the spirit of what Gubinelli, Ugurcan and Zachhuber did in \cite{GUZ}, and consider the modified paracontrolled expansion
	$$
	u=\PT_u^s X+u_s^\sharp
	$$
	with the truncated paraproduct $\PT^s$ defined below. For $s=s(\Xi)$ small enough, the map $\Phi^s(u):=u-\PT_u^sX$ is invertible as a perturbation of the identity and one can show that the subspace of such paracontrolled functions is indeed dense. The parameter $s$ will also be a very useful tool to investigate the different properties of $H$. Indeed, the Anderson operator will be given as
	$$
	Hu=Lu_s^\sharp+F_{\Xi,s}(u)
	$$
	with $F_{\Xi,s}:\CD(H)\subset L^2\to L^2$ an explicit operator and as $s$ goes to $0$, $u_s^\sharp$ gets closer to $u$ while $F_{\Xi,s}$ diverges. These different representations of $H$ will yield a family of bounds on the eigenvalues $\big(\lambda_n(\Xi)\big)_{n\ge1}$ of $H$ of the form
	$$
	m^-(\Xi,s)\lambda_n-m(\Xi,s)\le\lambda_n(\Xi)\le m^+(\Xi,s)\lambda_n+m(\Xi,s)
	$$
	with $(\lambda_n)_{n\ge1}$ the eigenvalues of $L$. In partiular, $m^-(\Xi,s)$ and $m^+(\Xi,s)$ converge to $1$ while $m(\Xi,s)$ diverges almost surely as $s$ goes to $0$. A particular choice for $s$ implies the simpler bounds
	$$
	\lambda_n-m_\delta^1(\Xi)\le\lambda_n(\Xi)\le(1+\delta)\lambda_n+m_\delta^2(\Xi)
	$$
	for any $\delta\in(0,1)$ which is expected but was not present in the initial work of Allez and Chouk \cite{AllezChouk} or the work of Labbé \cite{Labbe}. These bounds are sharp enough to imply the Weyl law
	$$
	\lim_{\lambda\to\infty}\lambda^{-1}\big|\{n\ge0;\lambda_n(\Xi)\le\lambda\}\big|=\frac{\textup{Vol}(M)}{4\pi}
	$$
	and give bounds for the tails of the eigenvalues.
	\item[$2)$] For the second question, one introduces the corrector $\DC$ with
	$$
	\PI(u,\xi)=u\PI(X,\xi)+\DC(u,X,\xi)+\PI(u^\sharp,\xi)
	$$
	for $u$ paracontrolled by $X$. One has to define the product $\PI(X,\xi)$ independently of the operator, this is the renormalisation step. To do so, we use the Wick product and set
	$$
	\PI(X,\xi):=\lim_{\varepsilon\to0}\Big(\PI(X_\varepsilon,\xi_\varepsilon)-\IE\big[\PI(X_\varepsilon,\xi_\varepsilon)\big]\Big)
	$$
	with $\xi_\varepsilon$ a regularisation of the noise. In some sense explained in Proposition \ref{Happrox}, the operator $H$ is the limit of the renormalised operators
	$$
	H_\varepsilon:=L+\xi_\varepsilon-c_\varepsilon
	$$
	with $c_\varepsilon:=\IE\big[\PI(X_\varepsilon,\xi_\varepsilon)\big]$ a smooth function diverging almost surely as $\varepsilon$ goes to $0$. Note that on the torus, the noise is invariant by translation and $c_\varepsilon$ is constant.
\end{enumerate}

The approach sketched above yields an operator $H:\CD(H)\subset L^2\to\CH^{2\alpha-2}$ with $\CD(H)$ the space of paracontrolled functions. In two dimensions, $2\alpha-2<0$ hence one needs to refine the definition of the domain to get an unbounded operator in $L^2$. To this purpose, Allez and Chouk introduced in \cite{AllezChouk} the subspace of $\CD(H)$ of strongly paracontrolled functions still dense in $L^2$. This was also used by Gubinelli, Ugurcan and Zachhuber in \cite{GUZ} and adapted to the dimension $3$ using a Hopf-Cole type transformation. We present here a different approach based on a higher order expansion. In particular, the domain of $H$ will consist of functions $u$ such that
$$
u=\PT_uX_1+\PT_uX_2+u^\sharp
$$
where $X_1\in\CC^\alpha,X_2\in\CC^{2\alpha}$ are noise-dependent functions and $u^\sharp\in\CH^2$. Note that since we want to get bounds in $\Xi$, quantitative estimates are needed and we keep track of the different explicit constants that appear, in particular how small $s$ needs to be with respect to the noise. If one is only interest in qualitative results, details of almost all computations can be skipped.

\smallskip

We shall first construct in Section $2.1$ the enhanced noise $\Xi$ from $\xi$ by a renormalisation procedure and prove exponential moments for its norm. The domain $\CD_\Xi$ of $H$ is constructed in Section $2.2$ and proved to be dense using a truncated paraproduct $\PT^s$. We show in particular in Proposition \ref{H2estimate} that the natural norms of $\CD_\Xi$ are equivalent to the norm operator; this will give the upper bound for the eigenvalues. Section $2.2$ is ended with the computation of the Hölder regularity of the elements of the domain. After showing that the operator is closed, we show in Section $2.3$ that $H$ is the limit of the operators $H_\varepsilon$ in some sense which yields the symmetry of $H$. We then control in Proposition \ref{H1estimate} the $\CH^1$ norm of $u^\sharp$ from the associated bilinear form applied to $u$; this will give the lower bound for the eigenvalues. This gives self-adjointness and pure point spectrum using the Babuška-Lax-Milgram theorem and we conclude the section with a bound on the convergence of the eigenvalues of $H_\varepsilon$ to $H$. Section $2.4$ treats the Schrödinger equation.

\medskip

As in the work of Allez and Chouk \cite{AllezChouk}, Labbé \cite{Labbe} and of Gubinelli, Ugurcan and Zachhuber \cite{GUZ}, we construct a dense random subspace of $L^2$ though a renormalisation step to get a self-adjoint operator with pure point spectrum. Our approach is different since we perform a second order expansion using paracontrolled calculus based on the heat semigroup on the manifold $M$. We refine the upper bounds on the eigenvalues obtained in \cite{AllezChouk} on the torus while also providing lower bounds. We get upper bounds for $\IP(\lambda_n(\Xi)\le\lambda)$ for $\lambda$ to $+\infty$ and $-\infty$. For $\lambda$ to $-\infty$, a bound was first given in \cite{Labbe} for a bounded domain with different boundary conditions. We have a more explicit dependence on $n$ while a less precise bound with respect to $\lambda$. To the best of our knowledge, no bounds for $\lambda$ to $+\infty$ were known. We also prove that the eigenfunctions of $H$ belong to $\CC^{1-}$ while the works \cite{AllezChouk,Labbe,GUZ} only gave Sobolev regularity. For the Schrödinger equation, we get on a manifold the same result as Gubinelli, Ugurcan and Zachhuber get on the torus, see \cite{GUZ}. As in their work, our construction of the Hamilton Anderson on $M$ could be used to study other evolution PDEs. All these results are new in our geometrical framework.

\bigskip

\subsection{Renormalisation}\label{subsecrenorm}

As explained in the introduction, an element of the domain of $H$ should behave like the linear part $X:=-L^{-1}\xi$ hence the product $u\xi$ does not make sense in two dimensions. Using the corrector, we are able to define the product $u\xi$ for $u$ paracontrolled by $X$ once the product $X\xi$ is defined. To do so, a naive approach would be to regularize the noise where $\xi_\varepsilon=\Phi(\varepsilon L)\xi$ is a regularisation of the noise and take $\varepsilon$ to $0$. The only condition we take is $\Phi$ such that $(\Phi(\varepsilon L))_\varepsilon$ belongs to the class $\mathsf{G}$, for example $\Phi(\varepsilon L)=e^{\varepsilon L}$ works. Since the product is ill-defined, the quantity $\PI(X_\varepsilon,\xi_\varepsilon)$ diverges as $\varepsilon$ goes to $0$ with $X_\varepsilon:=-L^{-1}\xi_\varepsilon$. The now usual way is to substract another diverging quantity $c_\varepsilon$ such that the limit
$$
\PI(X,\xi):=\lim_{\varepsilon\to0}\Big(\PI(X_\varepsilon,\xi_\varepsilon)-c_\varepsilon\Big)
$$
exists and take this as the definition of the product. This is the Wick renormalisation and the purpose of the following theorem with the renormalised Anderson Hamiltonian
$$
H_\varepsilon:=L+\xi_\varepsilon-c_\varepsilon.
$$
We insist here that $c_\varepsilon$ is not a priori a constant in the general geometrical setting of a two-dimensional manifold.

\medskip

\begin{theorem}
Let $\alpha<1$ and
$$
c_\varepsilon:=\IE\Big[\PI(X_\varepsilon,\xi_\varepsilon)\Big].
$$
Then there exists a random distribution $\PI(X,\xi)$ that belongs almost surely to $\CC^{2\alpha-2}$ and such that
$$
\lim_{\varepsilon\to0}\IE\Big[\big\|\PI(X,\xi)-(\PI(X_\varepsilon,\xi_\varepsilon)-c_\varepsilon)\big\|_{\CC^{2\alpha-2}}^p\Big]=0
$$
for any $p\ge1$.
\end{theorem}

\medskip

\begin{proof}
Since the noise is Gaussian, we only need to control second order moment using hypercontractivity. The resonant term $\PI(X_\varepsilon,\xi_\varepsilon)$ is a linear combination of terms of the form
$$
I_\varepsilon:=\int_0^1P_t^\bullet\left(Q_t^1X_\varepsilon\cdot Q_t^2\xi_\varepsilon\right)\frac{\drm t}{t}
$$
with $P\in\mathsf{StGC}^{[0,b]}$ and $Q^1,Q^2\in\mathsf{StGC}^{\frac{b}{2}}$. We also define the renormalised quantity
$$
J_\varepsilon:=I_\varepsilon-\IE[I_\varepsilon].
$$
Let $u\in(0,1),x\in M$ and $Q\in\mathsf{StGC}^r$ with $r>|2\alpha-2|$. The expectation $\IE\left[|Q_u\big(I_\varepsilon\big)(x)|^2\right]$ is given by the integral over $M^2\times[0,1]^2$ of
$$
K_{Q_uP_t^\bullet}(x,y)K_{Q_uP_s^\bullet}(x,z)\IE\Big[Q_t^1X_\varepsilon(y)Q_t^2\xi_\varepsilon(y)Q_s^1X_\varepsilon(z)Q_s^2\xi_\varepsilon(z)\Big]
$$
against the measure $\mu(\drm y)\mu(\drm z)(ts)^{-1}\drm t\drm s$. Using the Wick formula, we have
\begin{align*}
\IE&\Big[Q_t^1X_\varepsilon(y)Q_t^2\xi_\varepsilon(y)Q_s^1X_\varepsilon(z)Q_s^2\xi_\varepsilon(z)\Big]=\IE\left[Q_t^1X_\varepsilon(y)Q_t^2\xi_\varepsilon(y)\right]\IE\left[Q_s^1X_\varepsilon(z)Q_s^2\xi_\varepsilon(z)\right]\\
&\quad+\IE\left[Q_t^1X_\varepsilon(y)Q_s^1X_\varepsilon(z)\right]\IE\left[Q_t^2\xi_\varepsilon(y)Q_s^2\xi_\varepsilon(z)\right]+\IE\left[Q_t^1X_\varepsilon(y)Q_s^2\xi_\varepsilon(z)\right]\IE\left[Q_s^1X_\varepsilon(z)Q_t^2\xi_\varepsilon(y)\right]\\
&=(1)+(2)+(3)
\end{align*}
and this yields
$$
\IE\left[|Q_u\big(I_\varepsilon\big)(x)|^2\right]=I_\varepsilon^{(1)}(x)+I_\varepsilon^{(2)}(x)+I_\varepsilon^{(3)}(x).
$$
The first term corresponds exactly to the extracted diverging quantity since
$$
I_\varepsilon^{(1)}=\IE\left[\int_0^1Q_uP_t^\bullet\left(Q_t^1X_\varepsilon\cdot Q_t^2\xi_\varepsilon\right)\frac{\drm t}{t}\right]^2=\IE\big[Q_u(I_\varepsilon)\big]^2
$$
and we have
$$
\IE\left[|Q_u\big(J_\varepsilon\big)(x)|^2\right]=\IE\left[\Big\{Q_u\big(I_\varepsilon\big)(x)-\IE[Q_u\big(I_\varepsilon\big)](x)\Big\}^2\right]=I_\varepsilon^{(2)}(x)+I_\varepsilon^{(3)}(x).
$$
Using that $(\Psi(\varepsilon L))_\varepsilon$ belongs to $\mathsf{G}$, $\xi$ is an isometry from $L^2$ to square-integrable random variables and lemma \ref{lemmaPT}, we have
\begin{align*}
I_\varepsilon^{(2)}(x)+I_\varepsilon^{(3)}(x)&\lesssim\int_{M^2}\int_{[0,1]^2}K_{Q_uP_t^\bullet}(x,y)K_{Q_uP_s^\bullet}(x,z)\big\langle\CG_{2\varepsilon+t+s}(y,\cdot),\CG_{2\varepsilon+t+s}(z,\cdot)\big\rangle^2\mu(\drm y)\mu(\drm z)ts\drm t\drm s\\
&\lesssim\int_{M^2}\int_{[0,1]^2}K_{Q_uP_t^\bullet}(x,y)K_{Q_uP_s^\bullet}(x,z)\CG_{2\varepsilon+t+s}(y,z)^2\mu(\drm y)\mu(\drm z)ts\drm t\drm s\\
&\lesssim\int_{M^2}\int_{[0,1]^2}\CG_{u+t}(x,y)\CG_{u+s}(x,z)\CG_{2\varepsilon+t+s}(y,z)^2\mu(\drm y)\mu(\drm z)ts\drm t\drm s\\
&\lesssim\int_{M^2}\int_{[0,1]^2}(2\varepsilon+t+s)^{-\frac{d}{2}}\CG_{u+t}(x,y)\CG_{u+s}(x,z)\CG_{2\varepsilon+t+s}(y,z)\mu(\drm y)\mu(\drm z)ts\drm t\drm s\\
&\lesssim\int_{[0,1]^2}(2\varepsilon+t+s)^{-\frac{d}{2}}(\varepsilon+u+t+s)^{-\frac{d}{2}}ts\drm t\drm s\\
&\lesssim(\varepsilon+u)^{2-d}
\end{align*}
hence the family $\big(\PI(X_\varepsilon,\xi_\varepsilon)-c_\varepsilon\big)_{\varepsilon>0}$ is bounded in $\CC^{2\alpha-2}$ for any $\alpha<1$ since $d=2$. These computations also show that the associated linear combination of
$$
J:=\int_0^1\Big\{P_t^\bullet\left(Q_t^1X\cdot Q_t^2\xi\right)-\IE\left[P_t^\bullet\left(Q_t^1X\cdot Q_t^2\xi\right)\right]\Big\}\frac{\drm t}{t}
$$
yields a well-defined random distribution of $\CC^{2\alpha-2}$ for $\alpha<1$ that we denote $\PI(X,\xi)$. The same type of computations show the convergence and completes the proof.
\end{proof}

\medskip

The enhanced noise is defined as
$$
\Xi:=\big(\xi,\PI(X,\xi)\big)\in\CX^\alpha
$$ 
where $\CX^\alpha:=\CC^{\alpha-2}\times\CC^{2\alpha-2}$. One has to keep in mind that the notation $\PI(X,\xi)$ is only suggestive. In particular for almost every $\omega$, one has
$$
\PI\big(X,\xi\big)(\omega)\neq\PI\big(X(\omega),\xi(\omega)\big)
$$
since the product is almost surely ill-defined. We also denote the regularized enhanced noise $\Xi_\varepsilon:=\big(\xi_\varepsilon,\PI(X_\varepsilon,\xi_\varepsilon)-c_\varepsilon\big)$ with the norm
$$
\|\Xi-\Xi_\varepsilon\|_{\CX^\alpha}:=\|\xi-\xi_\varepsilon\|_{\CC^{\alpha-2}}+\big\|\PI(X,\xi)-\PI(X_\varepsilon,\xi_\varepsilon)+c_\varepsilon\big\|_{\CC^{2\alpha-2}}
$$
which goes to $0$ as $\varepsilon$ goes to $0$. Using that the noise is Gaussian and almost surely in $\CC^{-1-\kappa}$ for all $\kappa>0$, we have exponential moment for the norm of the enhanced noise.

\medskip

\begin{proposition}
There exists $h>0$ such that
$$
\IE\left[e^{h\|\xi\|_{\CC^{\alpha-2}}^2+h\|\PI(X,\xi)\|_{\CC^{2\alpha-2}}}\right]<\infty.
$$
\end{proposition}

\medskip

\begin{proof}
Let $t\in(0,1)$ and $Q\in\mathsf{StGC}^r$ with $r>|\alpha-2|$. Using the Gaussian hypercontractivity, we have
\begin{align*}
\IE\left[\|Q_t\xi\|_{L_x^p}^p\right]&=\int_M\IE\left[|Q_t\xi|^p(x)\right]\mu(\drm x)\\
&\le(p-1)^{\frac{p}{2}}\int_M\IE\left[|Q_t\xi|^2(x)\right]^{\frac{p}{2}}\mu(\drm x)
\end{align*}
hence we only need to bound the second moment, which is bounded by
$$
\IE\left[|Q_t\xi|^2(x)\right]=\|K_{Q_t}(x,\cdot)\|_{L^2}^2\lesssim\frac{1}{\mu\big(B(x,\sqrt{t})\big)}.
$$
Using that $\CB_{2p,2p}^{\alpha-2+\frac{1}{p}}\hookrightarrow\CB_{\infty,\infty}^{\alpha-2}$, we have
\begin{align*}
\IE\left[e^{h\|\xi\|_{\CC^{\alpha-2}}^2}\right]&=\sum_{p\ge0}\frac{h^p}{p!}\IE\big[\|\xi\|_{\CC^{\alpha-2}}^{2p}\big]\\
&\le\sum_{p=0}^{p_0}\frac{h^p}{p!}\IE\big[\|\xi\|_{\CC^{\alpha-2}}^{2p}\big]+\sum_{p>p_0}\frac{h^p}{p!}\IE\left[\|\xi\|_{\CB_{2p,2p}^{\alpha-2+\frac{1}{p}}}^{2p}\right]\\
&\lesssim\sum_{p=0}^{p_0}\frac{h^p}{p!}\IE\big[\|\xi\|_{\CC^{\alpha-2}}^{2p_0}\big]^{\frac{p}{p_0}}+\sum_{p>p_0}\frac{h^p(2p-1)^p}{p!}\textup{Vol}(M)
\end{align*}
for $p_0>\frac{2}{1-\alpha}$ hence the result for $h$ small enough. For the bound on $\PI(X,\xi)$, the computations are the same without the square since it belongs to the second Wiener chaos hence Gaussian hypercontractivity gives
$$
\IE\big[|Q_t\PI(X,\xi)|^p(x)\big]\le(p-1)^p\ \IE\big[|Q_t\PI(X,\xi)|^2(x)\big]^{\frac{p}{2}}.
$$
\end{proof}

\bigskip

\subsection{Domain of the Hamiltonian}

We first motivate the definition of the domain. Let $\alpha\in(\frac{2}{3},1)$ such that $\xi$ belongs almost surely to $\CC^{\alpha-2}$. Let $X\in\CC^\alpha$ be a noise-dependent function and consider $u=\PT_{u'}X+u^\sharp$ a function paracontrolled by $X$ with $u'\in\CH^\alpha$ and $u^\sharp\in\CH^{2\alpha}$. Then
\begin{align*}
Hu&=Lu+\xi u\\
&=L\big(\PT_{u'}X+u^\sharp\big)+\P_u\xi+\P_\xi u+\PI\big(\PT_{u'}X+u^\sharp,\xi\big)\\
&=\P_{u'}LX+\P_u\xi+\Big(Lu^\sharp+\P_\xi u+u'\PI(X,\xi)+\DC(u',X,\xi)+\PI(u^\sharp,\xi)\Big).
\end{align*}
Taking $u'=u$ and $-LX=\xi$, the first two terms cancel each other and we get
$$
Hu=Lu^\sharp+\P_\xi u+u\PI(X,\xi)+\DC(u,X,\xi)+\PI(u^\sharp,\xi)\in\CH^{2\alpha-2}.
$$
This yields an unbounded operator in $L^2$ with values in $\CH^{2\alpha-2}$. Since $2\alpha-2<0$, $Hu$ does not belong to $L^2$ hence we do not have an operator from $L^2$ to itself and this makes harder to study the spectral properties of $H$. To get around this, Allez and Chouk introduced in \cite{AllezChouk} the subspace of functions $u$ paracontrolled by $L^{-1}\xi$ such that $Hu$ does belong to $L^2$ called strongly paracontrolled functions. This approach was also used by Gubinelli, Ugurcan and Zachhuber in \cite{GUZ} however we proceed differently and use higher order expansions. Let $X_1:=X$ and $X_2\in\CC^{2\alpha}$ be another noise-dependent function. Given $u_2\in\CH^\alpha$ and $u^\sharp\in\CH^{3\alpha}$, we consider $u=\PT_uX_1+\PT_{u_2}X_2+u^\sharp$ and we have
\begin{align*}
Hu&=\P_{u_2}LX_2+u\PI(X_1,\xi)+\DC(u,X_1,\xi)+\P_{u_2}\PI(X_2,\xi)+\DD(u_2,X_2,\xi)\\
&\quad+\P_u\P_\xi X_1+\DS(u,X_1,\xi)+\P_\xi\PT_{u_2}X_2+\P_\xi u^\sharp+Lu^\sharp+\PI(u^\sharp,\xi).
\end{align*}
Taking $u_2=u$ and $-LX_2=\PI(X_1,\xi)+\P_\xi X_1$ cancels the terms of Sobolev regularity $2\alpha-2$ and we get
\begin{align*}
Hu&=\PI\big(u,\PI(X_1,\xi)\big)+\P_{\PI(X_1,\xi)}u+\DC(u,X_1,\xi)+\P_u\PI(X_2,\xi)+\DD(u,X_2,\xi)\\
&\quad+\DS(u,X_1,\xi)+\P_\xi\PT_uX_2+\P_\xi u^\sharp+Lu^\sharp+\PI(u^\sharp,\xi)
\end{align*}
hence $Hu\in\CH^{3\alpha-2}\subset L^2$. This motivates the following definition for the domain $\CD_\Xi$ of $H$ with 
$$
-LX_1:=\xi\quad\text{and}\quad-LX_2:=\PI(X_1,\xi)+\P_\xi X_1.
$$

\medskip

\begin{definition*}
We define the set $\CD_\Xi$ of functions paracontrolled by $\Xi$ as
$$
\CD_\Xi:=\big\{u\in L^2;\ u^\sharp:=u-\PT_uX_1-\PT_uX_2\in\CH^2\big\}.
$$
\end{definition*}

\medskip

The domain $\CD_\Xi$ is the random subspace of functions $u\in L^2$ paracontrolled by $X_1$ and $X_2$ up to a remainder $u^\sharp\in\CH^2$ given by the explicit formula
$$
u^\sharp=\Phi(u):=u-\PT_uX_1-\PT_uX_2.
$$
With this notation, we have $\CD_\Xi=\Phi^{-1}(\CH^2)$ and since $X_1+X_2\in\CC^\alpha$, we actually have $\CD_\Xi\subset\CH^\beta$ for every $\beta<\alpha$. However, we have no idea at this point if this domain is trivial or dense in $L^2$ and an inverse to $\Phi$ would be useful. However, it is not necessarily invertible so we introduce a parameter $s>0$ and consider
$$
\Phi^s:\left|\begin{array}{ccc}
\CD_\Xi&\to&\CH^2\\
u&\mapsto&u-\PT_u^s X_1-\PT_u^s X_2
\end{array}\right.
$$
where $\PT^s$ is defined as
$$
\PT_f^s g:=\sum_{\mathbf{a}\in\SA_b;a_2<\frac{b}{2}}\sum_{\mathbf{Q}\in\mathsf{StGC}^{\bf a}}b_{\bf Q}\int_0^s\widetilde Q_t^{1\bullet}\left(Q_t^2f\cdot\widetilde Q_t^3g\right)\frac{\drm t}{t}.
$$
The important property is that while still encoding the important information of the paraproduct $\PT$, the truncated paraproduct $\PT^s$ is small as an operator for $s$ small; this is quantified as follows and proved in Proposition \ref{ConvergeEta} in Appendix \ref{AppendixPC}.

\medskip

\begin{proposition}\label{BoundEta}
Let $\gamma\in(0,1)$ be a regularity exponent and $X\in\CC^\gamma$. For any $\beta\in[0,\gamma)$, we have
$$
\|u\mapsto\PT_u^s X\|_{L^2\to\CH^\beta}\lesssim\frac{s^{\frac{\gamma-\beta}{4}}}{\gamma-\beta}\|X\|_{\CC^\gamma}
$$
\end{proposition}

\medskip

Since $X_1$ and $X_2$ depends continuously on $\Xi$, this implies the existence of $m>0$ such that
$$
\|\PT_u^sX_1+\PT_u^sX_2\|_{\CH^\beta}\le m\frac{s^{\frac{\alpha-\beta}{4}}}{\alpha-\beta}\|\Xi\|_{\CX^\alpha}(1+\|\Xi\|_{\CX^\alpha})\|u\|_{L^2}
$$
thus the operator $u\mapsto\PT_u^s(X_1+X_2)$ is continuous from $L^2$ to $\CH^\beta$ for $\beta\in[0,\alpha)$ and arbitrary small as $s$ goes to $0$. Hence we get that
$$
\Phi^s:\CH^\beta\to\CH^\beta
$$ 
is invertible for $s=s(\Xi,\beta)$ small enough as a perturbation of the identity. Since $\PT_uX_i-\PT_u^sX_i$ is a smooth function for any $s>0$, the domain is still given by 
$$
\CD_\Xi=\Phi^{-1}(\CH^2)=(\Phi^s)^{-1}(\CH^2)
$$ 
and we have a decomposition given by $\Phi^s$ for any $u\in\CD_\Xi$, that is
$$
u=\PT_u^sX_1+\PT_u^sX_2+\Phi^s(u).
$$
In particular, we emphasize that \textbf{the domain does not depend on} $s$ while the decomposition we consider for element of the domain might. We denote
$$
\x:=\|\Xi\|_{\CX^\alpha}
$$
to keep track of the quantitative dependance with respect to the enhanced noise $\Xi$ and lighten the notation. We use the letter $x$ as a reminder of the noise dependance. For any $0\le\beta<\alpha$, we define
$$
s_\beta(\Xi):=\left(\frac{\alpha-\beta}{m\x(1+\x)}\right)^{\frac{4}{\alpha-\beta}}
$$
such that for $s<s_\beta(\Xi)$, the operator $\Phi^s:\CH^\beta\to\CH^\beta$ is invertible and we denote $\Gamma$ its inverse. We choose to drop the parameter $s$ in the notation to lighten the computations however the reader should keep in mind that the map $\Gamma$ depends on $s$. It is implicitly characterized by the relation
$$
\Gamma u^\sharp=\PT_{\Gamma u^\sharp}^sX_1+\PT_{\Gamma u^\sharp}^sX_2+u^\sharp
$$
for any $u^\sharp\in\CH^\beta$. Our choice of $\PT^s$ is motivated by the preservation of the intertwining relation
$$
\PT^s=L^{-1}\circ\P^s\circ L
$$
with $\P^s$ defined as $\PT^s$. The map $\Gamma$ will be a crucial tool to study the domain $\CD_\Xi$, in particular to show density in $L^2$. Continuity estimates for $\Phi^s$ and $\Gamma$ are given in the next proposition. Note that in the following, this bound of the form $\|a-b\|\le c$ will be used as $\|a\|\le\|b\|+c$ or $\|b\|\le\|a\|+c$.

\medskip

\begin{proposition}\label{GammaEstimates}
Let $\beta\in[0,\alpha)$ and $s\in(0,1)$. We have
$$
\|\Phi^s(u)-u\|_{\CH^\beta}\le\frac{m}{\alpha-\beta}s^{\frac{\alpha-\beta}{4}}\x(1+\x)\|u\|_{L^2}.
$$
If moreover $s<s_\beta(\Xi)$, this implies
$$
\|\Gamma u^\sharp\|_{\CH^\beta}\le\frac{1}{1-\frac{m}{\alpha-\beta}s^{\frac{\alpha-\beta}{4}}\x(1+\x)}\|u^\sharp\|_{\CH^\beta}.
$$
\end{proposition}

\medskip

\begin{proof}
The bounds on $\Phi^s$ follows directly from proposotion \ref{BoundEta}. Moreover since
$$
\frac{m}{\alpha-\beta}s^{\frac{\alpha-\beta}{4}}\x(1+\x)<1
$$
for $s<s_\beta(\Xi)$, the map $\Phi^s:\CH^\beta\to\CH^\beta$ is invertible and we have
$$
\|\Gamma u^\sharp\|_{\CH^\beta}\le\frac{1}{1-\frac{m}{\alpha-\beta}s^{\frac{\alpha-\beta}{4}}\x(1+\x)}\|u^\sharp\|_{\CH^\beta}.
$$
\end{proof}

\medskip

Let us insist that $\|u_s^\sharp\|_{\CH^\beta}$ is always controlled by $\|u\|_{\CH^\beta}$ while $s$ need to be small depending for $\|u\|_{\CH^\beta}$ to be controlled by $\|u_s^\sharp\|_{\CH^\beta}$. We also define the map $\Gamma_\varepsilon$ associated to the regularized noise $\Xi_\varepsilon$ as
$$
\Gamma_\varepsilon u^\sharp=\PT_{\Gamma_\varepsilon u^\sharp}^sX_1^{(\varepsilon)}+\PT_{\Gamma_\varepsilon u^\sharp}^sX_2^{(\varepsilon)}+u^\sharp
$$
with 
$$
-LX_1^{(\varepsilon)}:=\xi_\varepsilon\quad\text{and}\quad-LX_2^{(\varepsilon)}:=\PI(X_1^{(\varepsilon)},\xi_\varepsilon)-c_\varepsilon+\P_{\xi_\varepsilon}X_1^{(\varepsilon)}.
$$ 
It satisfies the same bound as $\Gamma$ with $\|\Xi_\varepsilon\|_{\CX^\alpha}$ and the following approximation lemma holds. We do not need to explicit the constant, it depends polynomialy on the noise $\Xi$ and diverges as $s$ goes to $s_\beta(\Xi)$.

\medskip

\begin{lemma}\label{GammaConvergence}
For any $0\le\beta<\alpha$ and $0<s<s_\beta(\Xi)$, we have
$$
\|\textup{Id}-\Gamma\Gamma_\varepsilon^{-1}\|_{L^2\to\CH^\beta}\lesssim_{\Xi,s,\beta}\|\Xi-\Xi_\varepsilon\|_{\CX^\alpha}.
$$
In particular, this implies the norm convergence of $\Gamma_\varepsilon$ to $\Gamma$ with the bound
$$
\|\Gamma-\Gamma_\varepsilon\|_{\CH^\beta\to\CH^\beta}\lesssim_{\Xi,s,\beta}\|\Xi-\Xi_\varepsilon\|_{\CX^\alpha}.
$$
\end{lemma}

\medskip

\begin{proof}
Given any $u\in\CH^\beta$, we have $u=\Gamma\Gamma^{-1}(u)=\Gamma(u-\PT_u^sX_1-\PT_u^sX_2)$. Using proposition \ref{GammaEstimates}, we get
\begin{align*}
\|u-\Gamma\Gamma_\varepsilon^{-1}(u)\|_{\CH^\beta}&=\big\|\Gamma\big(u-\PT_u^sX_1-\PT_u^sX_2\big)-\Gamma\big(u-\PT_u^sX_1^{(\varepsilon)}-\PT_u^sX_2^{(\varepsilon)}\big)\big\|_{\CH^\beta}\\
&=\Big\|\Gamma\Big(\PT_u^s\big(X_1^{(\varepsilon)}-X_1\big)+\PT_u^s\big(X_2^{(\varepsilon)}-X_2\big)\Big)\Big\|_{\CH^\beta}\\
&\le\frac{\alpha-\beta}{\alpha-\beta-ms^{\frac{\alpha-\beta}{4}}\x(1+\x)}\left\|\PT_u^s\big(X_1^{(\varepsilon)}-X_1\big)+\PT_u^s\big(X_2^{(\varepsilon)}-X_2\big)\right\|_{\CH^\beta}\\
&\lesssim\frac{s^{\frac{\alpha-\beta}{4}}(1+\x)}{\alpha-\beta-ms^{\frac{\alpha-\beta}{4}}\x(1+\x)}\|\Xi-\Xi_\varepsilon\|_{\CX^\alpha}\|u\|_{L^2}
\end{align*}
using the proposition \ref{BoundEta} and that $X_i^{(\varepsilon)}-X_i$ is $i$-linear in $\Xi_\varepsilon-\Xi$ for $i\in\{1,2\}$. The second statement follows from
$$
\|\Gamma_\varepsilon-\Gamma\|_{\CH^\beta\to\CH^\beta}=\|\left(\textup{Id}-\Gamma\Gamma_\varepsilon^{-1}\right)\Gamma_\varepsilon\|_{\CH^\beta\to\CH^\beta}\le\|\textup{Id}-\Gamma\Gamma_\varepsilon\|_{\CH^\beta\to\CH^\beta}\|\Gamma_\varepsilon\|_{\CH^\beta\to\CH^\beta}
$$
with the bound uniform in $\varepsilon$ for $s<s_\beta(\Xi_\varepsilon)$
$$
\|\Gamma_\varepsilon\|_{\CH^\beta\to\CH^\beta}\le\frac{\alpha-\beta}{\alpha-\beta-ms^{\frac{\alpha-\beta}{4}}\x(1+\x)}.
$$
\end{proof}

\medskip

This allows to prove density of the domain.

\medskip

\begin{corollary*}
The domain $\CD_\Xi$ is dense in $\CH^\beta$ for any $\beta\in[0,\alpha)$.
\end{corollary*}

\medskip

\begin{proof}
Given $f\in\CH^2$, $\Gamma(g_\varepsilon)\in\CD_\Xi$ where $g_\varepsilon=\Gamma_\varepsilon^{-1}f\in\CH^2$ thus we can conclude with the lemma \ref{GammaConvergence} that
$$
\lim_{\varepsilon\to0}\|f-\Gamma(g_\varepsilon)\|_{\CH^\beta}=0.
$$
The density of $\CH^2$ in $\CH^\beta$ then yields the result.
\end{proof}

\medskip

Taking into account in the previous computation the smooth term $e^{-L}$ coming from the intertwining relation, we are able to define $H$ as an unbounded operator in $L^2$ with domain $\CD_\Xi$ as follows.

\medskip

\begin{definition*}
We define the Anderson Hamiltonian $H:\CD_\Xi\to L^2$ as
$$
Hu=Lu^\sharp+\P_\xi u^\sharp+\PI(u^\sharp,\xi)+R(u)
$$
with $u^\sharp=\Phi(u)$ and $R:\CD_\Xi\to L^2$ given by
\begin{align*}
R(u)&:=\PI\big(u,\PI(X_1,\xi)\big)+\P_{\PI(X_1,\xi)}u+\DC(u,X_1,\xi)+\P_u\PI(X_2,\xi)+\DD(u,X_2,\xi)\\
&\quad+\DS(u,X_1,\xi)+\P_\xi\PT_uX_2-e^{-L}\left(\P_uX_1+\P_uX_2\right).
\end{align*}
\end{definition*}

\medskip

The parameter $s$ does not appear in the definition of $H$, it is a tool to study the properties of the operator. Indeed, one has different representations of $Hu$ as
$$
Hu=Lu_s^\sharp+\P_\xi u_s^\sharp+\PI(u_s^\sharp,\xi)+R(u)+\Psi^s(u)
$$
where $u_s^\sharp:=\Phi^s(u)$ and
$$
\Psi^s(u):=\Big(L+\P_\xi\cdot+\PI(\cdot,\xi)\Big)\big(\PT_u^s-\PT_u\big)(X_1+X_2).
$$
The different representations of $H$ through the parameter $s>0$ will be useful to get different bounds. For example, we can compare the graph norm of $H$ given as
$$
\|u\|_H^2:=\|u\|_{L^2}^2+\|Hu\|_{L^2}^2
$$
and the natural norms of the domain
$$
\|u\|_{\CD_\Xi}^2:=\|u\|_{L^2}^2+\|\Phi^s(u)\|_{\CH^2}^2
$$
with the following proposition. For $s\in(0,1)$ and $\delta>0$, we introduce the constant
$$
m_\delta^2(\Xi,s):=k\left(s^{\frac{\alpha-2}{2}}\x(1+\x^2)+s^{\frac{\alpha-\beta}{4}}\x^2(1+\x^{3})+\delta^{-3}\big(1+s^{\frac{\alpha}{4}}\x(1+\x)\big)\x^4(1+\x^8)\right)
$$
where the ``$2$'' refers to $\CH^2$ and for a constant $k>0$ large enough depending only on $M$ and $L$. In particular, it depends polynomialy on the enhanced noise and diverges as $s$ or $\delta$ goes to $0$.

\medskip

\begin{proposition}\label{H2estimate}
Let $u\in\CD_\Xi$ and $s>0$. For any $\delta>0$, we have
$$
(1-\delta)\|u_s^\sharp\|_{\CH^2}\le\|Hu\|_{L^2}+m_\delta^2(\Xi,s)\|u\|_{L^2}
$$
and
$$
\|Hu\|_{L^2}\le(1+\delta)\|u_s^\sharp\|_{\CH^2}+m_\delta^2(\Xi,s)\|u\|_{L^2}
$$
with $u_s^\sharp=\Phi^s(u)$.
\end{proposition}

\medskip

\begin{proof}
For any $s>0$, we have
$$
Hu=Lu_s^\sharp+\P_\xi u_s^\sharp+\PI(u_s^\sharp,\xi)+R(u)+\Psi^s(u).
$$
Then $Lu_s^\sharp\in L^2$ and for $\beta=\frac{1}{2}(\frac{2}{3}+\alpha)$, we have
\begin{align*}
\|R(u)\|_{L^2}&\lesssim\x(1+\x^2)\|u\|_{\CH^\beta}\\
\|\Psi^s(u)\|_{L^2}&\lesssim s^{\frac{\alpha-2}{2}}\x(1+\x^2)\|u\|_{L^2}\\
\|\P_\xi u_s^\sharp+\PI(u_s^\sharp,\xi)\|_{L^2}&\lesssim\|\xi\|_{\CC^{\alpha-2}}\|u_s^\sharp\|_{\CH^{\frac{4}{3}}}.
\end{align*}
One can bound the $\CH^\beta$ norm of $u$ using Proposition \ref{GammaEstimates} with
$$
\|u\|_{\CH^\beta}\le\|u_s^\sharp\|_{\CH^\beta}+\frac{m}{\alpha-\beta}s^{\frac{\alpha-\beta}{4}}\x(1+\x)\|u\|_{L^2}
$$
and since $\beta<1$, one has
$$
\|Lu_s^\sharp-Hu\|_{L^2}\lesssim\Big(s^{\frac{\alpha-2}{2}}\x(1+\x^2)+s^{\frac{\alpha-\beta}{4}}\x^2(1+\x^3)\Big)\|u\|_{L^2}+\x(1+\x^2)\|u_s^\sharp\|_{\CH^{\frac{4}{3}}}.
$$
Since $0<\beta<2$, we have for any $t>0$
\begin{align*}
\|u_s^\sharp\|_{\CH^{\frac{4}{3}}}&\lesssim\left\|\int_0^t (t'L)e^{-t'L}u_s^\sharp\frac{\drm t'}{t'}\right\|_{\CH^{\frac{4}{3}}}+\left\|e^{-tL}u_s^\sharp\right\|_{\CH^{\frac{4}{3}}}\\
&\lesssim t^{\frac{2}{3}}\|u_s^\sharp\|_{\CH^2}+t^{-\frac{4}{2}}\Big(1+s^{\frac{\alpha}{4}}x(1+x)\Big)\|u\|_{L^2}.
\end{align*}
Take
$$
t=\left(\frac{\delta}{k\x(1+\x^2)}\right)^{\frac{3}{2}}
$$
with $k$ the constant from the previous inequality and $\delta>0$. This yields 
$$
\|Lu_s^\sharp-Hu\|_{L^2}\lesssim m_\delta^2(\Xi,s)\|u\|_{L^2}+\delta\|u_s^\sharp\|_{\CH^2}.
$$
and completes the proof.
\end{proof}

\medskip

Finally, we can compute the Hölder regularity of the domain. In particular, this will implies the $\alpha$-Hölder regularity of the eigenfunctions of $H$.

\medskip

\begin{proposition*}
We have
$$
\CD_\Xi\subset\CC^{1-\kappa}
$$
for any $\kappa>0$.
\end{proposition*}

\medskip

\begin{proof}
The Besov embedding in two dimensions implies
$$
\CH^2\hookrightarrow\CB_{\infty,\infty}^1=\CC^1\hookrightarrow L^\infty
$$
and $\Phi^s:L^\infty\to L^\infty$ is also invertible hence
$$
\CD_\Xi=\big(\Phi^s\big)^{-1}(\CH^2)\subset L^\infty.
$$
Given any $u\in\CD_\Xi$, we get
\begin{align*}
\|u\|_{\CC^\alpha}&\lesssim\|u\|_{L^\infty}\|X_1+X_2\|_{\CC^\alpha}+\|u_s^\sharp\|_{\CC^\alpha}\\
&\lesssim_\Xi\|u\|_{L^\infty}+\|u_s^\sharp\|_{\CH^2}
\end{align*}
for any $\alpha<1$ and the proof is complete.
\end{proof}

\bigskip

\subsection{Self-adjointness and spectral properties}

We show that $H$ is a closed self-adjoint operator on its dense domain $\CD_\Xi\subset L^2$. This relies on approximation results and the Babuška-Lax-Milgram theorem. The spectrum is pure point and the eigenvalues verify a min-max principle that allows to get estimates depending on the eigenvalues of $L$. 

\medskip

\begin{proposition}
The operator $H$ is closed on its domain $\CD_\Xi$.
\end{proposition}

\medskip

\begin{proof}
Let $(u_n)_{n\ge0}\subset\CD_\Xi$ be a sequence such that 
$$
u_n\rightarrow u\quad\text{in }L^2\quad\text{and}\quad Hu_n\rightarrow v\quad\text{in }L^2.
$$
Proposition \ref{H2estimate} gives that $\big(\Phi(u_n)\big)_{n\ge0}$ is a Cauchy sequence in $\CH^2$ hence converges to $u^\sharp\in\CH^2$. Since $\Phi:L^2\to L^2$ is continuous, we have $\Phi(u)=u^\sharp$ hence $u\in\CD_\Xi$. Finally, we have
\begin{align*}
\|Hu-v\|_{L^2}&\le\|Hu-Hu_n\|_{L^2}+\|Hu_n-v\|_{L^2}\\
&\lesssim_\Xi\|u_n^\sharp-u^\sharp\|_{\CH^2}+\|u-u_n\|_{L^2}+\|Hu_n-v\|_{L^2}
\end{align*}
hence $Hu=v$ and $H$ is closed on $\CD_\Xi$.
\end{proof}

\medskip

In some sense, the operator $H$ should be the limit of the renormalised $H_\varepsilon$ as $\varepsilon$ goes to $0$. Since $\CD(H_\varepsilon)=\CH^2$, one can not compare directly the operators. However given any $u\in L^2$, we have
$$
u=\big(\Gamma\circ\Phi^s\big)(u)=\lim_{\varepsilon\to0}\big(\Gamma_\varepsilon\circ\Phi^s\big)(u)
$$
where we recall that $\Gamma_\varepsilon u^\sharp$ is implicitly defined for $u^\sharp\in L^2$ by the relation
$$ 
\Gamma_\varepsilon u^\sharp=\PT_{\Gamma_\varepsilon u^\sharp}^sX_1^{(\varepsilon)}+\PT_{\Gamma_\varepsilon u^\sharp}^sX_2^{(\varepsilon)}+u^\sharp.
$$
Thus for $u\in\CD_\Xi$, the approximation $u_\varepsilon:=\big(\Gamma_\varepsilon\circ\Phi^s\big)(u)
$ belongs to $\CH^2$ and one can consider the difference
$$
\|Hu-H_\varepsilon u_\varepsilon\|_{L^2}=\|(H\Gamma-H_\varepsilon\Gamma_\varepsilon)u^\sharp\|_{L^2}
$$
with $u^\sharp:=\Phi^s(u)$. The following proposition gives a bound for this quantity which yields the convergence as $\varepsilon$ goes to $0$ for $s$ is small enough. We do not need to explicit the constant, it depends polynomialy on the enhanced noise $\Xi$ and diverges as $s$ goes to $s_0(\Xi)$.

\medskip

\begin{proposition}\label{Happrox}
Let $u\in\CD_\Xi$ and $s>0$ small enough. Then
$$
\|Hu-H_\varepsilon u_\varepsilon\|_{L^2}\lesssim_{\Xi,s}\|u_s^\sharp\|_{\CH^2}\|\Xi-\Xi_\varepsilon\|_{\CX^\alpha}
$$
with $u_s^\sharp=\Phi^s(u)$ and $u_\varepsilon:=\Gamma_\varepsilon u_s^\sharp$. In particular, this implies that $H_\varepsilon\Gamma_\varepsilon$ converges to $H\Gamma$ in norm as $\varepsilon$ goes to $0$ as operators from $\CH^2$ to $L^2$.
\end{proposition}

\medskip

\begin{proof}
We have
$$
H_\varepsilon u_\varepsilon=Lu_s^\sharp+\P_{\xi_\varepsilon}u_s^\sharp+\PI(u_s^\sharp,\xi)+R_\varepsilon(u_\varepsilon)+\Psi_\varepsilon^s(u_\varepsilon)
$$
where $R_\varepsilon$ and $\Psi_\varepsilon^s$ are defined as $R$ and $\Psi^s$ with $\Xi_\varepsilon$ instead of $\Xi$. For $\beta=\frac{1}{2}(\frac{2}{3}+\alpha)$, we have
\begin{align*}
\|R(u)-&R_\varepsilon(u_\varepsilon)\|_{L^2}\le\|R(u-u_\varepsilon)\|_{L^2}+\|(R-R_\varepsilon)(u_\varepsilon)\|_{L^2}\\
&\lesssim \x(1+\x^2)\|u-u_\varepsilon\|_{\CH^\beta}+(1+\x)\|\Xi-\Xi_\varepsilon\|_{\CX^\alpha}\|u_\varepsilon\|_{\CH^\beta}\\
&\lesssim\Big(\x(1+\x^2)\|\Gamma-\Gamma_\varepsilon\|_{\CH^\beta\to\CH^\beta}+(1+\x)\|\Gamma_\varepsilon\|_{\CH^\beta\to\CH^\beta}\|\Xi-\Xi_\varepsilon\|_{\CX^\alpha}\Big)\|u_s^\sharp\|_{\CH^2}
\end{align*}
and the same reasoning gives
$$
\|\Psi^s(u)-\Psi_\varepsilon^s(u)\|_{L^2}\lesssim_{s,\Xi}\|\Xi-\Xi_\varepsilon\|_{\CX^\alpha}.
$$
Thus one completes the proof with the bound $\|\Gamma-\Gamma_\varepsilon\|_{\CH^\beta\to\CH^\beta}$ from Lemma \ref{GammaConvergence}.
\end{proof}

\medskip

The symmetry of $H$ immediately follows.

\medskip

\begin{corollary*}
The operator $H$ is symmetric.
\end{corollary*}

\medskip

\begin{proof}
Let $u,v\in\CD_\Xi$ and consider $u^\sharp:=\Phi^s(u)$ and $v^\sharp:=\Phi^s(v)$ for $s<s_0(\Xi)$. Since $H_\varepsilon$ is a symmetric operator, we have
$$
\langle Hu,v\rangle=\lim_{\varepsilon\to0}\langle H_\varepsilon\Gamma_\varepsilon u^\sharp,\Gamma_\varepsilon v^\sharp\rangle=\lim_{\varepsilon\to0}\langle\Gamma_\varepsilon u^\sharp,H_\varepsilon\Gamma_\varepsilon v^\sharp\rangle=\langle u,Hv\rangle
$$
using that $H_\varepsilon\Gamma_\varepsilon$ converges to $H\Gamma$ and $\Gamma_\varepsilon$ to $\Gamma$ in norm convergence.
\end{proof}

\medskip

The next proposition states that the quadratic form associated to $H$ is bounded from below by the $\CH^1$ norm of $u^\sharp$. This weak coercivity property will give below self-adjointness with the Babuška-Lax-Milgram theorem. This was already used in the work \cite{GUZ} of Gubinelli, Ugurcan and Zachhuuber, where the proof of self-adjointness relies on the reasoning of almost duality encoded in the operator $\DO$. For $s\in(0,1)$ and $\delta>0$, introduce the constant
\begin{align*}
m_\delta^1(\Xi,s)&:=k\Big\{x(1+\x^2)+s^{\frac{\alpha-\beta}{4}}\x^2(1+\x^3)+s^{\frac{\alpha-2}{2}}\x(1+\x^2)+s^{\frac{\alpha-4}{2}}\x\\
&\quad+\delta^{-\frac{\beta}{1-\beta}}\Big(\x(1+\x^2)+s^{\frac{\alpha-\beta}{4}}\x^2(1+\x)\Big)^{\frac{\beta}{1-\beta}}\Big(1+s^{\frac{\alpha}{4}}\x(1+\x)\Big)\Big\}
\end{align*}
where $\beta=\frac{1}{2}(\frac{2}{3}+\alpha)$ and for a constant $k>0$ large enough depending only on $M$ and $L$ while the ``$1$'' refers to $\CH^1$. In particular, it depends polynomialy on the enhanced noise and diverges as $s$ or $\delta$ goes to $0$.

\medskip

\begin{proposition}\label{H1estimate}
Let $u\in\CD_\Xi$ and $s>0$. For any $\delta>0$, we have
$$
(1-\delta)\langle\nabla u_s^\sharp,\nabla u_s^\sharp\rangle\le\langle u,Hu\rangle+m_\delta^1(\Xi,s)\|u\|_{L^2}^2
$$
and
$$
(1-\delta)\langle\nabla u_s^\sharp,\nabla u_s^\sharp\rangle\le\langle u,H_\varepsilon u\rangle+m_\delta^1(\Xi,s)\|u\|_{L^2}^2
$$
where $u_s^\sharp=\Phi^s(u)$.
\end{proposition}

\medskip

\begin{proof}
For $u\in\CD_\Xi$, we have
$$
Hu=Lu_s^\sharp+\P_\xi u_s^\sharp+\PI(u_s^\sharp,\xi)+R(u)+\Psi^s(u)
$$
with $u_s^\sharp=\Phi^s(u)\in\CH^2$. Thus
\begin{align*}
\langle u,Lu_s^\sharp\rangle&=\big\langle\PT_u^sX_1,Lu_s^\sharp\big\rangle+\big\langle\PT_u^sX_2,Lu_s^\sharp\big\rangle+\big\langle u_s^\sharp,Lu_s^\sharp\big\rangle\\
&=\big\langle\P_u^s LX_1,u_s^\sharp\big\rangle+\big\langle\P_u^s LX_2,u_s^\sharp\big\rangle+\big\langle\nabla u_s^\sharp,\nabla u_s^\sharp\big\rangle
\end{align*}
and this yields
\begin{align*}
\langle u,&Hu\rangle=-\big\langle\P_u^s\xi,u_s^\sharp\big\rangle+\big\langle\P_u^s LX_2,u_s^\sharp\big\rangle+\big\langle\nabla u_s^\sharp,\nabla u_s^\sharp\big\rangle+\big\langle u,\P_\xi u_s^\sharp+\PI(u_s^\sharp,\xi)\big\rangle+\big\langle u,R(u)+\Psi^s(u)\big\rangle\\
&=-\DO(u,\xi,u_s^\sharp)+\big\langle\P_u^s LX_2,u_s^\sharp\big\rangle+\big\langle\nabla u_s^\sharp,\nabla u_s^\sharp\big\rangle+\big\langle u,\P_\xi u_s^\sharp\big\rangle+\big\langle u,R(u)+\Psi^s(u)\big\rangle+\big\langle(\P_u-\P_u^s)\xi,u_s^\sharp\big\rangle
\end{align*}
where $\DO(u,\xi,u^\sharp)=\big\langle\P_u\xi,u^\sharp\big\rangle-\big\langle u,\PI(u^\sharp,\xi)\big\rangle$. For $\beta:=\frac{1}{2}(\frac{2}{3}+\alpha)$, we have
\begin{align*}
\big|\big\langle u,R(u)\big\rangle\big|&\lesssim\|u\|_{L^2}\|R(u)\|_{L^2}\lesssim\x(1+\x^2)\|u\|_{L^2}\|u\|_{\CH^\beta},\\
\big|\big\langle u,\P_\xi u_s^\sharp\big\rangle\big|&\lesssim\|u\|_{\CH^\beta}\|\P_\xi u_s^\sharp\|_{\CC^{2\beta-2}}\lesssim\x\|u\|_{\CH^\beta}\|u_s^\sharp\|_{\CH^\beta},\\
\big|\big\langle\P_uLX_2,u_s^\sharp\big\rangle\big|&\lesssim\|\P_uLX_2\|_{\CH^{2\beta-2}}\|u_s^\sharp\|_{\CH^\beta}\lesssim\x^2\|u\|_{L^2}\|u_s^\sharp\|_{\CH^\beta}.
\end{align*}
Using Proposition \ref{Duality}, we have
$$
\big|\DO(u,\xi,u_s^\sharp)\big|\lesssim\|\xi\|_{\CC^{\alpha-2}}\|u\|_{\CH^\beta}\|u_s^\sharp\|_{\CH^\beta}\lesssim\x\|u\|_{\CH^\beta}\|u_s^\sharp\|_{\CH^\beta}.
$$
Finally, we have
\begin{align*}
\big|\big\langle u,\Psi^s(u)\big\rangle\big|&\lesssim\|u\|_{L^2}\|\Psi^s(u)\|_{L^2}\lesssim s^{\frac{\alpha-2}{2}}\x(1+\x^2)\|u\|_{L^2}^2\\
\big|\big\langle(\P_u-\P_u^s)\xi,u_s^\sharp\big\rangle\big|&\lesssim\|(\P_u-\P_u^s)\xi\|_{L^2}\|u_s^\sharp\|_{L^2}\lesssim s^{\frac{\alpha-4}{2}}\x\|u\|_{L^2}\|u_s^\sharp\|_{L^2}
\end{align*}
with Proposition \ref{DivergeEta} in Appendix \ref{AppendixPC}. Since $u\in\CD_\Xi$, we have
$$
\|u\|_{\CH^\beta}\le\|u_s^\sharp\|_{\CH^\beta}+\frac{m}{\alpha-\beta}s^{\frac{\alpha-\beta}{4}}\x(1+\x)\|u\|_{L^2}
$$
hence there exists $k>0$ such that
\begin{align*}
\big\langle\nabla u_s^\sharp,\nabla u_s^\sharp\big\rangle\le&\big\langle u,Hu\big\rangle+k\Big(\x(1+\x^2)+s^{\frac{\alpha-\beta}{4}}\x^2(1+\x^3)+s^{\frac{\alpha-2}{2}}\x(1+\x^2)+s^{\frac{\alpha-4}{2}}\x\Big)\|u\|_{L^2}^2\\
&+k\Big(\x(1+\x^2)+s^{\frac{\alpha-\beta}{4}}\x^2(1+\x)\Big)\|u_s^\sharp\|_{\CH^\beta}.
\end{align*}
Since $0<\beta<1$, we have for any $t>0$
\begin{align*}
\|u_s^\sharp\|_{\CH^\beta}^2&\lesssim\left\|\int_0^t (t'L)e^{-t'L}u_s^\sharp\frac{\drm t'}{t'}\right\|_{\CH^\beta}^2+\left\|e^{-tL}u_s^\sharp\right\|_{\CH^\beta}^2\\
&\lesssim t^{1-\beta}\|u_s^\sharp\|_{\CH^1}^2+t^{-\beta}\Big(1+s^{\frac{\alpha}{4}}\x(1+\x)\Big)^2\|u\|_{L^2}^2.
\end{align*}
Given any $\delta>0$, we set 
$$
t=\left(\frac{\delta}{k'\Big(\x(1+\x^2)+s^{\frac{\alpha-\beta}{4}}\x^2(1+\x)\Big)}\right)^{\frac{1}{1-\beta}}
$$ 
where $k'>0$ the constant from the previous inequality and this yields
$$
(1-\delta)\big\langle\nabla u_s^\sharp,\nabla u_s^\sharp\big\rangle\le\big\langle u,Hu\big\rangle+m_\delta^1(\Xi,s)\|u\|_{L^2}.
$$
The same computations show
$$
(1-\delta)\langle\nabla u_s^\sharp,\nabla u_s^\sharp\rangle\le\langle u,H_\varepsilon u\rangle+m_\delta^1(\Xi_\varepsilon,s)\|u\|_{L^2}^2.
$$
Since $\|\Xi_\varepsilon-\Xi\|_\alpha$ goes to $0$ as $\varepsilon$ goes to $0$, the result holds uniformly in $\varepsilon$ with $m_\delta^1(\Xi,s)$.
\end{proof}

\medskip

This implies that $H$ is almost surely bounded below by the random variable $-m_\delta^1(\Xi,s)$ for any $\delta>0$ and $s>0$. Using the Babuška-Lax-Milgram theorem, one gets an invertible operator via the solution of
$$
(H+k_\Xi)u=v
$$
for $k_\Xi>m_\delta^1(\Xi,s)$ and $v\in L^2$.

\medskip

\begin{proposition}\label{Resolvent}
Let $\delta\in(0,1)$ and $s>0$. Then for any constant $k_\Xi>m_\delta^1(\Xi,s)$, the operators $H+k_\Xi$ and $H_\varepsilon+k_\Xi$ are invertible. Moreover the operators
\begin{align*}
\big(H+k_\Xi\big)^{-1}&:L^2\to\CD_\Xi\\
\big(H_\varepsilon+k_\Xi\big)^{-1}&:L^2\to\CH^2
\end{align*}
are bounded.
\end{proposition}

\medskip

\begin{proof}
We want to use the theorem of Babuška-Lax-Milgram, see \cite{Babuska}. This is a generalization of the Lax-Milgram theorem with a weaker condition of coercivity. Since $k_\Xi>m_\delta^1(\Xi,s)$, Proposition \ref{H1estimate} gives
$$
\big(k_\Xi-m_\delta^1(\Xi,s)\big)\|u\|_{L^2}^2<\big\langle(H+k_\Xi)u,u\big\rangle
$$
for $u\in\CD_\Xi$. Considering the norm
$$
\|u\|_{\CD_\Xi}^2=\|u\|_{L^2}^2+\|u_s^\sharp\|_{\CH^2}^2
$$
on $\CD_\Xi$, this yields a weakly coercive operator using Proposition \ref{H2estimate} in the sense that
$$
\|u\|_{\CD_\Xi}\lesssim_\Xi\|(H+k_\Xi)u\|_{L^2}=\sup_{\|v\|_{L^2}=1}\big\langle(H+k_\Xi)u,v\big\rangle
$$
for any $u\in\CD_\Xi$. Moreover, the bilinear map
$$
\left.\begin{array}{cccc}
B:&\CD_\Xi\times L^2&\to&\IR\\
&(u,v)&\mapsto&\big\langle(H+k_\Xi)u,v\big\rangle	
\end{array}\right.
$$
is continuous since Proposition \ref{H2estimate} implies
$$
\quad|B(u,v)|\le\|(H+k_\Xi)u\|_{L^2}\|v\|_{L^2}\lesssim_\Xi\|u\|_{\CD_\Xi}\|v\|_{L^2}
$$
for $u\in\CD_\Xi$ and $v\in L^2$. The last condition we need is that for any $v\in L^2\backslash\{0\}$, we have
$$
\sup_{\|u\|_{\CD_\Xi}=1}|B(u,v)|>0.
$$
Let assume that there exists $v\in L^2$ such that $B(u,v)=0$ for all $u\in\CD_\Xi$. Then
$$
\forall u\in\CD_\Xi,\quad\langle u,v\rangle_{\CD_\Xi,\CD_\Xi^*}=0.
$$
hence $v=0$ as an element of $\CD_\Xi^*$. By density of $\CD_\Xi$ in $L^2$, this implies $v=0$ in $L^2$ hence the property we want. By the theorem of Babuška-Lax-Milgram, for any $f\in L^2$ there exists a unique $u\in\CD_\Xi$ such that
$$
\forall v\in L^2,\quad B(u,v)=\langle f,v\rangle.
$$
Moreover, we have $\|u\|_{\CD_\Xi}\lesssim_\Xi\|f\|_{L^2}$ hence the result for $(H+k_\Xi)^{-1}$. The same argument works for $H_\varepsilon+k_\Xi$ since proposition \ref{H1estimate} also holds for $H_\varepsilon$ with bounds uniform in $\varepsilon$.
\end{proof}

\medskip

Using that a closed symmetric operator on a Hilbert space is self-adjoint if it has at least one real value in its resolvent set, this immediatly implies that $H$ and $H_\varepsilon$ are self-adjoint, see \cite{RS}. Moreover, the resolvant is a compact operator from $L^2$ to itself since $\CD_\Xi\subset\CH^\beta$ for any $\beta\in[0,\alpha)$ hence the following result.

\medskip

\begin{corollary}\label{SpectralResult}
The operators $H$ and $H_\varepsilon$ are self-adjoint with discret spectrum $\big(\lambda_n(\Xi)\big)_{n\ge1}$ and $\big(\lambda_n(\Xi_\varepsilon)\big)_{n\ge1}$ which are nondecreasing diverging sequences without accumulation points. Moreover, we have
$$
L^2=\underset{n\ge1}{\bigoplus}\ \textup{Ker}\big(H-\lambda_n(\Xi)\big)
$$
with each kernel being of finite dimension. We finally have the min-max principle
$$
\lambda_n(\Xi)=\inf_D\sup_{u\in D;\|u\|_{L^2}=1}\langle Hu,u\rangle
$$
where $D$ is any $n$-dimensional subspace of $\CD_\Xi$ that can also be given as
$$
\lambda_n(\Xi)=\sup_{v_1,\ldots,v_{n-1}\in L^2}\ \inf_{\underset{\|u\|_{L^2}=1}{u\in\textup{Vect}(v_1,\ldots,v_{n-1})^\bot}}\langle Hu,u\rangle.
$$
\end{corollary}

\medskip

A natural question now is to estimate the size of the eigenvalues of $H$ and try to get back geometric informations on the manifold $M$ as one can do from the Laplacian. Let $\lambda$ be an eigenvalue of $H$ and $u\in\CD_\Xi$ such that
$$
Hu=\lambda u.
$$
Then there exists $u^\sharp\in\CH^2$ such that $u=\Gamma u^\sharp$ thus
$$
H\Gamma u^\sharp=\lambda\Gamma u^\sharp.
$$
This yields
$$
H\Gamma u^\sharp=\lambda u^\sharp+\lambda\big(\Gamma-\textup{Id}\big)u^\sharp
$$
hence one can relate the spectrum of $H$ to the one of $H\Gamma$ and the parameter $s$ measures the error since
$$
\big(\Gamma-\textup{Id}\big)u^\sharp=\PT_{\Gamma u^\sharp}^sX_1+\PT_{\Gamma u^\sharp}^sX_2.
$$
And since $H\Gamma$ is a perturbation of $L$, one can relate the spectrum of $H\Gamma$ to the spectrum of $L$, as stated in the following proposition using the min-max result. We denote by $(\lambda_n)_{n\ge1}$ the non-decreasing positive sequence of the eigenvalues of $L$, since it corresponds to the case $\Xi=0$. For $s\in(0,1)$ and $\delta>0$, introduce the constant
$$
m_\delta^+(\Xi,s):=(1+\delta)\Big(1+\frac{m}{\alpha}s^{\frac{\alpha}{4}}\x(1+\x)\Big).
$$
If $s<s_0(\Xi)$, we also introduce
$$
m_\delta^-(\Xi,s):=(1-\delta)\frac{1}{1-\frac{m}{\alpha}s^{\frac{\alpha}{4}}\x(1+\x)}.
$$
In particular, the constants depend polynomialy on the enhanced noise $\Xi$ and converge to $1$ as $\delta$ and $s$ goes to $0$. Moreover, $m_\delta^-(\Xi,s)$ diverges as $s$ goes to $s_0(\Xi)$. Write $a,b\le c$ to mean that we have both $a\le c$ and $b\le c$.

\medskip

\begin{proposition}\label{EVestimates}
Let $s\in(0,1)$ and $\delta>0$. Given any $n\in\IZ^+$, we have
$$
\lambda_n(\Xi),\lambda_n(\Xi_\varepsilon)\le m_\delta^+(\Xi,s)\lambda_n+1+\frac{m}{\alpha}s^{\frac{\alpha}{4}}\x(1+\x)+m_\delta^2(\Xi,s).
$$
If moreover $s<s_0(\Xi)$, we have
$$
\lambda_n(\Xi),\lambda_n(\Xi_\varepsilon)\ge m_\delta^-(\Xi,s)\lambda_n-m_\delta^1(\Xi,s).
$$
\end{proposition}

\medskip

\begin{proof}
Let $u_1^\sharp,\ldots,u_n^\sharp\in\CH^2$ be an orthonormal family of eigenfunctions of $L$ associated to $\lambda_1,\ldots,\lambda_n$ and consider
$$
u_i:=\Gamma u_i^\sharp\in\CD_\Xi
$$
for $1\le i\le n$. Since $\Gamma$ is invertible, the family $(u_1,\ldots,u_n)$ is free thus the min-max representation of $\lambda_n(\Xi)$ yields
$$
\lambda_n(\Xi)\le\sup_{\underset{\|u\|_{L^2}=1}{u\in\textup{Vect}(u_1,\ldots,u_n)}}\langle Hu,u\rangle.
$$
Given any normalised $u\in\textup{Vect}(u_1,\ldots,u_n)$, we have
$$
\langle Hu,u\rangle\le\|Hu\|_{L^2}\le(1+\delta)\|u_s^\sharp\|_{\CH^2}+m_\delta^2(\Xi,s)
$$
for $u_s^\sharp=\Phi^s(u)$ using Proposition \ref{H2estimate}. Moreover
$$
\|u_s^\sharp\|_{\CH^2}\le(1+\lambda_n)\|u_s^\sharp\|_{L^2}\le(1+\lambda_n)\Big(1+\frac{m}{\alpha}s^{\frac{\alpha}{4}}\x(1+\x)\Big)
$$
hence the upper bound
$$
\lambda_n(\Xi)\le m_\delta^+(\Xi,s)\lambda_n+1+\frac{m}{\alpha}s^{\frac{\alpha}{4}}\x(1+\x)+m_\delta^2(\Xi,s).
$$
For the lower bound, we use the min-max representation of $\lambda_n(\Xi)$ under the form
$$
\lambda_n(\Xi)=\sup_{v_1,\ldots,v_{n-1}\in L^2}\ \inf_{\underset{\|u\|_{L^2}=1}{u\in\textup{Vect}(v_1,\ldots,v_{n-1})^\bot}}\langle Hu,u\rangle.
$$
Introducing
$$
F:=\textup{Vect}(u_m;m\ge n),
$$ 
we have that $F^\bot$ is a subspace of $L^2$ of finite dimension $n-1$ thus there exists a orthogonal family $(v_1,\ldots,v_{n-1})$ such that $F^\bot=\textup{Vect}(v_1,\ldots,v_{n-1})$. Since $F$ is a closed subspace of $L^2$ as an intersection of hyperplans, we have $F=\textup{Vect}(v_1,\ldots,v_{n-1})^\bot$ hence
$$
\lambda_n(\Xi)\ge\inf_{\underset{\|u\|_{L^2}=1}{u\in F}}\langle Hu,u\rangle.
$$
Let $u\in F$ with $\|u\|_{L^2}=1$. Using Proposition \ref{H1estimate}, we have
\begin{align*}
\langle Hu,u\rangle&\ge(1-\delta)\langle\nabla u_s^\sharp,\nabla u_s^\sharp\rangle-m_\delta^1(\Xi,s)\\
&\ge(1-\delta)\langle u_s^\sharp,Lu_s^\sharp\rangle-m_\delta^1(\Xi,s)\\
&\ge(1-\delta)\lambda_n\|u_s^\sharp\|_{L^2}^2-m_\delta^1(\Xi,s).
\end{align*}
Finally using Proposition \ref{GammaEstimates} for $s<s_0(\Xi)$, we get
$$
\langle Hu,u\rangle\ge\frac{1-\delta}{1-\frac{m}{\alpha}s^{\frac{\alpha}{4}}\x(1+\x)}\lambda_n-m_\delta^1(\Xi,s)
$$
and the proof is complete.
\end{proof}

\medskip

There is a wide range of choices for the constants $s\in(0,1)$ and $\delta>0$. For example, one can take
$$
s=\left(\frac{\alpha\delta}{m\x(1+\x)}\right)^{\frac{4}{\alpha}}
$$
for any $\delta\in(0,1)$ and get
$$
\lambda_n-m_\delta^1(\Xi)\le\lambda_n(\Xi)\le(1+\delta)\lambda_n+m_\delta^2(\Xi)
$$
for explicit constants $m_\delta^1$ and $m_\delta^2$, where the lower bound holds since $\delta<1$ gives $s<s_0(\Xi)$. This implies the following estimate for the tail of all the eigenvalues. A more precise result of this type was already obtained in \cite{Labbe} by Labbé in the flat case for $\lambda$ to $-\infty$ with $a=1$ where he also obtained a lower bound on the convergence of the form
$$
e^{-a_n\lambda}\le\IP(\lambda_n(\Xi)\le -\lambda)\le e^{-b_n\lambda}
$$
for $\lambda>0$ large enough and $a_n>b_n>0$ two constants. Here we get upper bounds for $\lambda$ to $+\infty$ and $-\infty$.

\medskip

\begin{corollary}\label{TailEstimates}
For any $n\in\IZ^+$ and $\lambda\in\IR$, we have
$$
1-me^{-h(\lambda-2\lambda_n)^{\frac{1}{12}}}\le\IP\big(\lambda_n(\Xi)\le\lambda\big)\le me^{-h(\lambda_n-\lambda)^{\frac{1}{5}}}
$$
where $m=\IE\big[e^{h\|\Xi\|_{\CX^\alpha}}\big]$.
\end{corollary}

\medskip

\begin{proof}
Fix $\delta\in(0,1)$ and let $\lambda\in\IR$. Denote $m_1=m_\delta^1$ and $m_2=m_\delta^2$. We have
$$
\IP\big(\lambda_n(\Xi)\le\lambda\big)\le\IP\big(\lambda_n-m_1\le\lambda\big)
$$
and
$$
\IP\big(\lambda_n(\Xi)>\lambda\big)\le\IP\big((1+\delta)\lambda_n+m_2>\lambda\big)
$$
thus
$$
\IP\big(m_2\le\lambda-(1+\delta)\lambda_n\big)\le\IP\big(\lambda_n(\Xi)\le\lambda\big)\le\IP\big(m_1\ge-\lambda+\lambda_n\big).
$$
There exists two constants $a_1,a_2>0$ such that
$$
m_i\le 1+\|\Xi\|_{\CX^\alpha}^{a_i}
$$
for $i\in\{1,2\}$, take for example $a_1=5$ and $a_2=12$. Hence
\begin{align*}
\IP\big(m_i\ge y\big)&=\IP\big(\|\Xi\|_{\CX^\alpha}\ge(y-1)^{\frac{1}{a_i}}\big)\\
&=\IP\big(e^{h\|\Xi\|_{\CX^\alpha}}\ge e^{hy^{\frac{1}{a_i}}}\big)\\
&\le e^{-hy^{\frac{1}{a_i}}}\IE\big[e^{h\|\Xi\|_{\CX^\alpha}}\big]
\end{align*}
using Markov inequality and this yields
$$
1-me^{-h(\lambda-(1+\delta)\lambda_n)^{\frac{1}{a_2}}}\le\IP\big(\lambda_n(\Xi)\le\lambda\big)\le me^{-h(\lambda_n-\lambda)^{\frac{1}{a_1}}}
$$
where $m=\IE\big[e^{h\|\Xi\|_{\CX^\alpha}}\big]$.
\end{proof}

\medskip

We proved that $H_\varepsilon$ converges to $H$ is some sense as $\varepsilon$ goes to $0$. The following proposition gives the convergence of $H_\varepsilon+k_\Xi$ to $H+k_\Xi$ in resolvent sense as $\varepsilon$ goes to $0$. We do not need to explicit the constant, it depends polynomialy on the enhanced noise $\Xi$.

\medskip

\begin{proposition}\label{ResolventConvergence}
Let $s\in(0,1)$ and $\delta>0$. Then for any constant $k_\Xi>m_\delta^1(\Xi,s)$ and $\beta\in[0,\alpha)$, we have
$$
\|(H_\varepsilon+k_\Xi)^{-1}-(H+k_\Xi)^{-1}\|_{L^2\to\CH^\beta}\lesssim_{\Xi,\beta}\|\Xi-\Xi_\varepsilon\|_{\CX^\alpha}.
$$
In particular, $(H_\varepsilon+k_\Xi)^{-1}$ converges to $(H+k_\Xi)^{-1}$ in norm as operator from $L^2$ to itself.
\end{proposition}

\medskip

\begin{proof}
Let $v\in L^2$. Since $H+k_\Xi:\CD_\Xi\to L^2$ is invertible, there exists $u\in\CD_\Xi$ such that 
$$
v=(H+k_\Xi)u
$$ 
thus
$$
\|(H+k_\Xi)^{-1}v-(H_\varepsilon+k_\Xi)^{-1}v\|_{L^2}=\|u-(H_\varepsilon+k_\Xi)^{-1}(H+k_\Xi)u\|_{L^2}.
$$
We introduce $u_\varepsilon:=\Gamma_\varepsilon\Phi^s(u)$ which converges to $u$ in $L^2$ and we have
$$
\|u-(H_\varepsilon+k_\Xi)^{-1}(H+k_\Xi)u\|_{L^2}\le\|u-u_\varepsilon\|_{L^2}+\|u_\varepsilon-(H_\varepsilon+k_\Xi)^{-1}(H+k_\Xi)u\|_{L^2}.
$$
Since Lemma \ref{GammaConvergence} gives
$$
\|u-u_\varepsilon\|_{L^2}\lesssim_{\Xi,s}\|\Xi-\Xi_\varepsilon\|_{\CX^\alpha},
$$
we only have to bound the second term. We have
\begin{align*}
\|u_\varepsilon-(H_\varepsilon+k_\Xi)^{-1}(H+k_\Xi)u\|_{L^2}&=\|(H_\varepsilon+k_\Xi)^{-1}\big((H_\varepsilon+k_\Xi)u_\varepsilon-(H+k_\Xi)u\big)\|_{L^2}\\
&\lesssim\|(H_\varepsilon+k_\Xi)u_\varepsilon-(H+k_\Xi)u\|_{L^2}\\
&\lesssim\|H_\varepsilon u_\varepsilon-Hu\|_{L^2}+k\|u_\varepsilon-u\|_{L^2}
\end{align*}
using Proposition \ref{Happrox}. In the end, we have
$$
\|(H+k)^{-1}v-(H_\varepsilon+k)^{-1}v\|_{L^2}\lesssim\|u_s^\sharp\|_{\CH^2}\|\Xi-\Xi_\varepsilon\|_{\CX^\alpha}
$$
hence the result since $(H+k)^{-1}:L^2\to\CD_\Xi$ is continuous.
\end{proof}

\medskip

This allows to get a bound on the convergence of $\lambda_n(\Xi_\varepsilon)$ to $\lambda_n(\Xi)$ as $\varepsilon$ goes to $0$.

\medskip

\begin{corollary}
For all $n\in\IN^*$, we have
$$
\left|\frac{1}{\lambda_n(\Xi)+k_\Xi}-\frac{1}{\lambda_n(\Xi_\varepsilon)+k_\Xi}\right|\lesssim_\Xi\|\Xi-\Xi_\varepsilon\|_{\CX^\alpha}.
$$
In particular, this implies
$$
|\lambda_n(\Xi)-\lambda_n(\Xi_\varepsilon)|\lesssim_\Xi(\lambda_n(\Xi)+k_\Xi)^2\ \|\Xi-\Xi_\varepsilon\|_{\CX^\alpha}
$$
\end{corollary}

\medskip

\begin{proof}
We use the min-max principle for $(H+k_\Xi)^{-1}$ and $(H_\varepsilon+k_\Xi)^{-1}$ and denote $\mu_n$ and $\mu_n^{(\varepsilon)}$ their $n$-th smallest eigeinvalue with multiplicity. Let $D_n=\textup{Vect}(v_1,\ldots,v_n)$ with $v_i$ an eigenfunction associated to $\mu_i^{(\varepsilon)}$ for $1\le i\le n$. Then for all $v\in D_n$ with $\|v\|_{L^2}=1$, we have
\begin{align*}
\big\langle(H+k_\Xi)^{-1}u,u\big\rangle&=\Big\langle\big((H+k_\Xi)^{-1}-(H_\varepsilon+k_\Xi)^{-1}\big)u,u\Big\rangle+\big\langle(H_\varepsilon+k_\Xi)^{-1}u,u\big\rangle\\
&\le\big\|(H+k_\Xi)^{-1}-(H_\varepsilon+k_\Xi)^{-1}\big\|_{L^2\to L^2}+\mu_n^{(\varepsilon)}
\end{align*}
hence with proposition \ref{ResolventConvergence} we get
$$
\mu_n-\mu_n^{(\varepsilon)}\lesssim_\Xi\|\Xi-\Xi_\varepsilon\|_{\CX^\alpha}.
$$
Using the same argument with eigeinfunctions associated to $(H+k_\Xi)^{-1}$, we get
$$
|\mu_n-\mu_n^{(\varepsilon)}|\lesssim_\Xi\|\Xi-\Xi_\varepsilon\|_{\CX^\alpha}.
$$
Thus this gives
$$
\left|\frac{1}{\lambda_n(\Xi)+k_\Xi}-\frac{1}{\lambda_n(\Xi_\varepsilon)+k_\Xi}\right|\lesssim_\Xi\|\Xi-\Xi_\varepsilon\|_{\CX^\alpha}
$$
and completes the proof with the upper bound on $\lambda_n(\Xi)$.
\end{proof}

\medskip

We conclude this section by giving as corollary the Weyl law for the Anderson Hamiltonian $H$.

\medskip

\begin{corollary}\label{weyl}
We have
$$
\lim_{\lambda\to\infty}\lambda^{-1}|\{n\ge0;\lambda_n\le\lambda\}|=\frac{\textup{Vol}(M)}{4\pi}.
$$
\end{corollary}

\medskip

\begin{proof}
Proposition \ref{EVestimates} gives the bounds
$$
\lambda_n-m_\delta^1(\Xi)\le\lambda_n(\Xi)\le(1+\delta)\lambda_n+m_\delta^2(\Xi)
$$
for any $\delta\in(0,1)$. Using the lower bound, $\lambda_n(\Xi)\le\lambda$ implies
$$
\lambda_n\le\lambda+m_\delta^1(\Xi,s)
$$
and thus
$$
|\{n\ge1;\lambda_n(\Xi)\le\lambda\}|\le N\big(\lambda+m_\delta^1(\Xi,s)\big)
$$
where $N(x)$ denotes the number of eigenvalues of $L$ smaller than $x$. Using the upper bound, we get
$$
N\left(\frac{\lambda-m_\delta^2(\Xi)}{1+\delta}\right)\le\big|\{n\ge1;\lambda_n(\Xi)\le\lambda\}\big|\le N\big(\lambda+m_\delta^1(\Xi)\big)
$$
hence the proof is complete using the result for the Laplace-Beltrami operator.
\end{proof}

\bigskip

\subsection{Stochastic nonlinear Schrödinger equation}\label{subsecSCH}

The construction of the Anderson Hamiltonian allows the study of associated evolution equations. This was the motivation for the work \cite{GUZ} of Gubinelli, Ugurcan and Zachhuber and they studied the nonlinear Schrödinger and wave equations on the torus in two and three dimensions, see the references therein for other approaches. Our work allows to do the same on a two-dimensional manifold. As an example, we give results for the cubic nonlinear Schrödinger equation associated to $H$. See the work \cite{DebusscheWeber} of Debussche and Weber for the equation on the torus where they use a Hopf-Cole type transformation. This was extended in \cite{TV} by Tzvetkov and Visciglia to the fourth order nonlinearity. 

\smallskip

Define the positive operator
$$
H^+:=H+k_\Xi
$$
with $k_\Xi$ as in Proposition \ref{H1estimate}. Proposition \ref{Resolvent} yields a characterization of the domain and the form domain which is defined as follows.

\medskip

\begin{definition*}
We define the form domain of $H$ denoted $\CD_\Xi(\sqrt{H^+})$ as the closure of the domain under the norm
$$
\|u\|_{\CD_\Xi(\sqrt{H^+})}:=\sqrt{\langle u,H^+u\rangle}
$$
\end{definition*}

\medskip

\begin{proposition}\label{EquivNorm}
For $s<s_0(\Xi)$ and $u\in L^2$,
$$
\Big(u\in\CD_\Xi(H^+)\Big)\iff\Big(u_s^\sharp=\Phi^s(u)\in\CH^2\Big).
$$ 
and we have the bounds
$$
\|u_s^\sharp\|_{\CH^2}\lesssim_{\Xi,s}\|H^+u\|_{L^2}\lesssim_{\Xi,s}\|u_s^\sharp\|_{\CH^2}.
$$
Moreover, we have
$$
\Big(u\in\CD_\Xi(\sqrt{H^+})\Big)\iff\Big(\Phi^s(u)=u_s^\sharp\in\CH^1\Big)
$$ 
with the bounds
$$
\|u_s^\sharp\|_{\CH^1}\lesssim_{\Xi,s}\|u\|_{\CD_\Xi(\sqrt{H^+})}\lesssim_{\Xi,s}\|u_s^\sharp\|_{\CH^1}.
$$
\end{proposition}

\medskip

\begin{proof}
As stated, Proposition \ref{H1estimate} yields
$$
\|u_s^\sharp\|_{\CH^1}\lesssim_{\Xi,s}\|u\|_{\CD_\Xi(\sqrt{H^+})}.
$$
In fact, the inequality that is proved is
$$
\big|\langle H^+u,u\rangle-\langle\nabla u_s^\sharp,\nabla u_s^\sharp\rangle\big|\le k_\Xi\|u\|_{L^2}+\delta\|u_s^\sharp\|_{\CH^1}
$$
thus one also get the other estimate
$$
\|u\|_{\CD_\Xi(\sqrt{H^+})}\lesssim_{\Xi,s}\|u_s^\sharp\|_{\CH^1}.
$$
\end{proof}

\medskip

This yields a version of Brezis-Gallouët inequality for the Anderson Hamiltonian. In some sense, it interpolates the $L^\infty$-norm between the energy norm and the logarithm of the domain norm. This was already obtained in \cite{GUZ} by Gubinelli, Ugurcan and Zachhuber on the torus.

\medskip

\begin{theorem}
For any $v\in\CD_\Xi(H^+)$, we have
$$
\|v\|_{L^\infty}\lesssim_\Xi\|v\|_{\CD_\Xi(\sqrt{H^+})}\left(1+\sqrt{\log\left(1+\frac{\|v\|_{\CD_\Xi(H^+)}}{\|v\|_{\CD(\sqrt{H^+})}}\right)}\right).
$$
For any $v\in\CH^2$, we have
$$
\|v\|_{L^\infty}\lesssim_\Xi\|\sqrt{H^+_\varepsilon}v\|_{L^2}\left(1+\sqrt{\log\left(1+\frac{\|H^+_\varepsilon v\|_{L^2}}{\|\sqrt{H^+_\varepsilon}v\|_{L^2}}\right)}\right).
$$
In particular, the second inequality holds uniformly in $\varepsilon$.
\end{theorem}

\medskip

\begin{proof}
For any $t>0$, we have
$$
\|v\|_{L^\infty}\le\left\|\int_0^t(t'L)e^{-t'L}v\frac{\drm t'}{t'}\right\|_{L^\infty}+\|e^{-tL}v\|_{L^\infty}.
$$
One can bound the integral as
\begin{align*}
\left\|\int_0^t(t'L)e^{-t'L}v\frac{\drm t'}{t'}\right\|_{L^\infty}&\lesssim\int_0^t\|Lv\|_{L^2}\drm t'\\
&\lesssim t\|v\|_{\CH^2}
\end{align*}
and the remainder as
\begin{align*}
\|e^{-tL}v\|_{L^\infty}&\lesssim\left\|\int_t^1(t'L)e^{-t'L}v\frac{\drm t'}{t'}\right\|_{L^\infty}+\|e^{-L}v\|_{L^\infty}\\
&\lesssim\left(\int_t^1\frac{\drm t'}{t'}\right)^{\frac{1}{2}}\left(\int_t^1\|(t'L)e^{-t'L}v\|_{L^\infty}^2\frac{\drm t'}{t'}\right)^{\frac{1}{2}}+\|v\|_{\CH^1}\\
&\lesssim\left(\int_t^1\frac{\drm t'}{t'}\right)^{\frac{1}{2}}\left(\int_t^1(t')^{-1}\|(t'L)e^{-t'L}v\|_{L^2}^2\frac{\drm t'}{t'}\right)^{\frac{1}{2}}+\|v\|_{\CH^1}\\
&\lesssim\|v\|_{\CH^1}\big(1+|\log(t)|^{\frac{1}{2}}\big),
\end{align*}
to get
$$
\|v\|_{L^\infty}\lesssim t\|v\|_{\CH^2}+\big(1+|\log(t)|^{\frac{1}{2}}\big)\|v\|_{\CH^1}.
$$
Taking $\|v\|_{\CH^1}\le 1$ and $t=\frac{\sqrt{\log(1+\|v\|_{\CH_2})}}{1+\|v\|_{\CH^2}}>0$, we get the classical Brezis-Gallouet inequality, that is
$$
\|v\|_{L^\infty}\lesssim1+\sqrt{\log\left(1+\|v\|_{\CH^2}\right)}.
$$
Thus for $\|v\|_{\CD(\sqrt{H^+})}\le1$, we have
\begin{align*}
\|v\|_{L^\infty}&\lesssim_\Xi\|v^\sharp\|_{L^\infty}\\
&\lesssim_\Xi1+\sqrt{\log\left(1+\|v^\sharp\|_{\CH^2}\right)}\\
&\lesssim_\Xi1+\sqrt{\log\left(1+\|H^+v\|_{\CD(H^+)}\right)}
\end{align*}
using proposition \ref{EquivNorm}. Since every estimates also hold for $H^+_\varepsilon$ with bound uniform in $\varepsilon$, we also get the estimate for $H^+_\varepsilon$. Applying this result to $\frac{v}{\|v\|_{\CD(\sqrt{H^+})}}$ yields the general inequality.
\end{proof}

\medskip

This inequality can be used for example to study the cubic nonlinear Schrödinger equation with multiplicative noise
$$
i\partial_tu+Lu+u\xi=-|u|^2u
$$
with initial condition $u_0\in\CD_\Xi$. The construction of the operator $H$ immediatly yields the renormalised solution $u(t,\cdot):=e^{-itH}u_0$ to the linear equation
$$
i\partial_tu+Lu+u\xi=0
$$
given any $u_0\in\CD_\Xi$ as done in \cite{GUZ}. Remark that when one regularizes the question, one also has to consider a suitable sequence of initial data $(u_0^{(\varepsilon)})_{\varepsilon>0}$, it is often refered to as ``well-prepared data'' in the litterature. This can also be used to solve the associated equation with cubic nonlinearity even if we can not apply the same theorem as Brezis and Gallouët in \cite{BG} since we do not have a control on the cubic term from $\CD_\Xi$ to itself. One could modify the domain taking into account the term $\PI(X_1,X_1)$ in $X_2$ to get a domain stable by multiplication. However since a direct computation as done by Gubinelli, Ugurcan and Zachhuber in \cite{GUZ} is enough, it is not necessary. In particular, the proof of the following theorem works exactly as in their work and is left to the reader.

\medskip

\begin{theorem}
Let $T>0$ and $u_0\in\CD_\Xi$. Then there exists a unique solution $u\in C\big([0,T],\CD(T)\big)\cap C^1\big([0,T],L^2\big)$ to the equation
$$
\left\{\begin{array}{ccc}
i\partial_tu&=&H^+u-|u|^2u\\
u(0,\cdot)&=&u_0
\end{array}\right.\quad\text{on }[0,T]\times M.
$$
Moreover, $u$ is the $L^2$-limit of the solutions $u_\varepsilon\in C\big([0,T],\CH^2\big)\cap C^1\big([0,T],L^2\big)$ of solutions to the equations
$$
\left\{\begin{array}{ccc}
i\partial_tu_\varepsilon&=&H_\varepsilon^+ u_\varepsilon-|u_\varepsilon|^2u_\varepsilon\\
u_\varepsilon(0,\cdot)&=&u_0^{(\varepsilon)}
\end{array}\right.\quad\text{on }[0,\infty[\times M,
$$
with the initial data 
$$
u_0^{(\varepsilon)}:=(H_\varepsilon^+)^{-1}H^+u_0\in\CH^2
$$
which converges to $u_0$ in $L^2$. We also have the convergences
\begin{align*}
u_\varepsilon(t)&\to u(t)\quad\text{in }L^2,\\
H_\varepsilon^+u_\varepsilon(t)&\to H^+u(t)\quad\text{in }L^2,\\
\partial_tu_\varepsilon(t)&\to\partial_tu(t)\quad\text{in }L^2
\end{align*}
for all $t\in[0,T]$.
\end{theorem}

\medskip

\begin{remark}
From the solution to
$$
i\partial_tu=H^+u-|u|^2u,
$$
on the torus, one easily gets the solution to the initial equation
$$
i\partial_tv=Hv-|v|^2v
$$
via the change of variable $u(t,\cdot)=e^{tk_\Xi}v(t,\cdot)$ since $k_\Xi$ is a constant. One could want to do the same in a manifold setting and compare the initial regularized equation
$$
i\partial_tu=Lu+\xi_\varepsilon u-|u|^2u.
$$
with the renormalised equation
$$
i\partial_tv=Lv+\xi_\varepsilon v-c_\varepsilon v-|v|^2v
$$
as Tzvetkov and Visciglia's Theorem $1.1$ from \cite{TV}. It is not clear what the change of variable should be on a manifold since $c_\varepsilon$ is a function and not a constant. It should still be possible to find an appropriate change of variable even though this requires some work.
\end{remark}

\bigskip

\appendix

\section{Approximation operators} \label{AppendixApprox}

We describe in this Appendix technical estimates needed in our continuous setting analog of the discrete Paley-Littlewood decomposition. The following proposition is the analog of the inclusions of $\ell^p$ spaces.

\medskip

\begin{proposition}
Let $p,q_1,q_2\in[1,\infty]$ with $q_1\le q_2$. For $f\in L^p$ and $\alpha\in\IR$, we have
$$
\left\|t^{-\frac{\alpha}{2}}\|Q_tf\|_{L_x^p}\right\|_{L^{q_2}(t^{-1}\drm t)}\lesssim\left\|t^{-\frac{\alpha}{2}}\|Q_tf\|_{L_x^p}\right\|_{L^{q_1}(t^{-1}\drm t)}.
$$
\end{proposition}

\medskip

\begin{proof}
We prove $\|\cdot\|_{L^\infty(t^{-1}\drm t)}\lesssim\|\cdot\|_{L^q(t^{-1}\drm t)}$ for any $q\in[1,\infty)$ and the result follows from duality. To get this, we use
$$
Q_t=2\int_{\frac{t}{2}}^tQ_s\left(\frac{t}{s}\right)^{a+1}P_{t-s}^{(c)}\frac{\drm s}{s}
$$
for any $Q\in\mathsf{StGC}^a$ and $t\in(0,1]$ which yields
$$
\|Q_tf\|_{L^p}\lesssim\int_{\frac{t}{2}}^t\|Q_sf\|_{L^p}\frac{\drm t}{t}\lesssim\left(\int_{\frac{t}{2}}^t\|Q_sf\|_{L^p}^q\frac{\drm t}{t}\right)^{\frac{1}{q}}.
$$
\end{proof}

\medskip

One needs the following bound to keep an accurate track of the constant in different estimates.

\medskip

\begin{lemma}
Let $r>0$ and $\alpha\in(-r,r)$. We have
$$
\int_0^\infty\left(\frac{u}{1+u^2}\right)^ru^\alpha\frac{\drm u}{u}\le\frac{2r}{r^2-\alpha^2}.
$$
\end{lemma}

\medskip

\begin{proof}
Since
$$
1=\frac{1+u^2}{1+u^2}=\frac{1}{1+u^2}+\frac{u^2}{1+u^2}
$$
and $u\ge0$, we have
\begin{align*}
\int_0^\infty\left(\frac{u}{1+u^2}\right)^ru^\alpha\frac{\drm u}{u}&=\int_0^1\left(\frac{u}{1+u^2}\right)^ru^\alpha\frac{\drm u}{u}+\int_1^\infty\left(\frac{u}{1+u^2}\right)^ru^\alpha\frac{\drm u}{u}\\
&\le\frac{1}{r+\alpha}+\frac{1}{r-\alpha}
\end{align*}
hence the bound.
\end{proof}

\medskip

The next lemma describes the localisation of the cancellation in our continuous context, including the dependance on $s>0$.

\medskip

\begin{lemma}
Let $r>0$ and $\alpha\in(-r,r)$. Given any $q\in[1,\infty]$, we have
$$
\left\|u^{-\alpha}\int_0^1\left(\frac{tu}{(t+u)^2}\right)^rf(t)\frac{\drm t}{t}\right\|_{L^q(u^{-1}\drm u)}\le\frac{2r}{r^2-\alpha^2}\left\|u^{-\alpha}f(u)\right\|_{L^q(u^{-1}\drm u)}.
$$
We also have
$$
\left\|u^{-\alpha}\int_0^s\left(\frac{tu}{(t+u)^2}\right)^rf(t)\frac{\drm t}{t}\right\|_{L^q(u^{-1}\drm u)}\le\frac{2r}{r^2-\alpha^2}s^{\beta-\alpha}\left\|u^{-\beta}f(u)\right\|_{L^q(u^{-1}\drm u)}
$$
for any $s>0$ and $\beta\in(\alpha,r)$.
\end{lemma}

\medskip

\begin{proof}
For $q=\infty$, we have
\begin{align*}
\left|\int_0^1\left(\frac{tu}{(t+u)^2}\right)^rf(t)\frac{\drm t}{t}\right|&\le\|t^{-\alpha}f(t)\|_{L^\infty}\int_0^1\left(\frac{tu}{(t+u)^2}\right)^rt^\alpha\frac{\drm t}{t}\\
&\le\left(\int_0^\infty\left(\frac{v}{1+v^2}\right)^rv^\alpha\frac{\drm v}{v}\right)u^\alpha\|t^{-\alpha}f(t)\|_{L^\infty}\\
&\le\frac{2r}{r^2-\alpha^2}u^\alpha\|t^{-\alpha}f(t)\|_{L^\infty}
\end{align*}
which yields the result. For $q=1$, we have
\begin{align*}
\int_0^1u^{-\alpha}\left|\int_0^1\left(\frac{tu}{(t+u)^2}\right)^rf(t)\frac{\drm t}{t}\right|\frac{\drm u}{u}&\le\int_0^1\left(\int_0^1\left(\frac{tu}{(t+u)^2}\right)^ru^{-\alpha}\frac{\drm u}{u}\right)|f(t)|\frac{\drm t}{t}\\
&\le\left(\int_0^\infty\left(\frac{v}{1+v^2}\right)^rv^\alpha\frac{\drm v}{v}\right)\int_0^1t^{-\alpha}|f(t)|\frac{\drm t}{t}\\
&\le\frac{2r}{r^2-\alpha^2}\int_0^1t^{-\alpha}|f(t)|\frac{\drm t}{t}.
\end{align*}
The result then follows for any $q\in(1,\infty)$ by interpolation. For the dependance with respect to $s$, we also interpolate between $q=1$ and $q=\infty$ and conclud with
\begin{align*}
\left|\int_0^s\left(\frac{tu}{(t+u)^2}\right)^rf(t)\frac{\drm t}{t}\right|&\le\|t^{-\beta}f(t)\|_{L^\infty}\int_0^s\left(\frac{tu}{(t+u)^2}\right)^rt^\beta\frac{\drm t}{t}\\
&\le s^{\beta-\alpha}\|t^{-\beta}f(t)\|_{L^\infty}\int_0^s\left(\frac{tu}{(t+u)^2}\right)^rt^\alpha\frac{\drm t}{t}\\
&\le\frac{2r}{r^2-\alpha^2}s^{\beta-\alpha}u^\alpha\|t^{-\alpha}f(t)\|_{L^\infty}
\end{align*}
and
\begin{align*}
\int_0^1u^{-\alpha}\left|\int_0^s\left(\frac{tu}{(t+u)^2}\right)^rf(t)\frac{\drm t}{t}\right|\frac{\drm u}{u}&\le\int_0^s\left(\int_0^1\left(\frac{tu}{(t+u)^2}\right)^ru^{-\alpha}\frac{\drm u}{u}\right)|f(t)|\frac{\drm t}{t}\\
&\le\frac{2r}{r^2-\alpha^2}\int_0^st^{-\alpha}|f(t)|\frac{\drm t}{t}\\
&\le\frac{2r}{r^2-\alpha^2}s^{\beta-\alpha}\int_0^1 t^{-\beta}|f(t)|\frac{\drm t}{t}.
\end{align*}
\end{proof}

\medskip

Finally, we have the following estimate for integrals.

\medskip

\begin{lemma}
Given any $\alpha>0$ and $q\in[1,\infty]$, we have
$$
\left\|u^{-\frac{\alpha}{2}}\int_0^uf(t)\frac{\drm t}{t}\right\|_{L^q(u^{-1}\drm u)}\le\frac{2}{\alpha}\|u^{-\frac{\alpha}{2}}f(u)\|_{L^q(u^{-1}\drm u)}.
$$
\end{lemma}

\medskip

\begin{proof}
We proceed again by interpolation proving the estimate for $q=\infty$ and $q=1$. Using that $\alpha>0$, we have
$$
\left|\int_0^uf(t)\frac{\drm t}{t}\right|\le\|t^{-\frac{\alpha}{2}}f(t)\|_{L^\infty}\int_0^ut^{\frac{\alpha}{2}}\frac{\drm t}{t}\le\frac{2}{\alpha}u^{\frac{\alpha}{2}}\|t^{-\frac{\alpha}{2}}f(t)\|_{L^\infty}
$$
and
$$
\int_0^1u^{-\frac{\alpha}{2}}\left|\int_0^uf(t)\frac{\drm t}{t}\right|\frac{\drm u}{u}\le\int_0^1\left(\int_t^1u^{-\frac{\alpha}{2}}\frac{\drm u}{u}\right)|f(t)|\frac{\drm t}{t}\le\frac{2}{\alpha}\int_0^1 t^{-\frac{\alpha}{2}}|f(t)|\frac{\drm t}{t}.
$$
\end{proof}

\bigskip

\section{Paracontrolled calculus} \label{AppendixPC}

We give in this Appendix proofs of estimates needed in paracontrolled calculus. We shall first prove the estimates for the paraproduct $\P$ and resonant operator $\PI$ in Sobolev spaces. It works as for Hölder spaces with $L^2$ estimates instead of $L^\infty$.

\medskip

\begin{proposition}\label{ParaSobolevAppendix}
Let $\alpha,\beta\in(-2b,2b)$ be regularity exponent.
\begin{itemize}
	\item[$\centerdot$] If $\alpha>0$, then $(f,g)\mapsto\P_fg$ is continuous from $\CC^\alpha\times\CH^\beta$ to $\CH^\beta$ and from $\CH^\alpha\times\CC^\beta$ to $\CH^\beta$.
	\item[$\centerdot$] If $\alpha<0$, then $(f,g)\mapsto\P_fg$ is continuous from $\CC^\alpha\times\CH^\beta$ to $\CH^{\alpha+\beta}$ and from $\CH^\alpha\times\CC^\beta$ to $\CH^{\alpha+\beta}$.
	\item[$\centerdot$] If $\alpha+\beta>0$, then $(f,g)\mapsto\PI(f,g)$ is continuous from $\CH^\alpha\times\CC^\beta$ to $\CH^{\alpha+\beta}$.
\end{itemize}
\end{proposition}

\medskip

\begin{proof}
Let $f\in\CH^\alpha$ and $g\in\CC^\beta$ with $\alpha<0$. We want to compute the regularity $\CH^{\alpha+\beta}$ of $\P_fg$ hence let $Q\in\mathsf{StGC}^r$ with $r>|\alpha+\beta|$. Recall that $\P_fg$ is a linear combination of terms of the form
$$
\int_0^1Q_t^{1\bullet}\left(P_tf\cdot Q_t^2g\right)\frac{\drm t}{t}
$$
with $Q^1,Q^2\in\mathsf{StGC}^{\frac{b}{2}}$ and $P\in\mathsf{StGC}^{[0,b]}$. Given $s\in(0,1]$, we have
\begin{align*}
\left\|\int_0^1 Q_sQ_t^{1\bullet}\left(P_tf\cdot Q_t^2g\right)\frac{\drm t}{t}\right\|_{L_x^2}&\lesssim\int_0^1\left(\frac{ts}{(t+s)^2}\right)^{\frac{r}{2}}\left\|P_tf\cdot Q_t^2g\right\|_{L_x^2}\frac{\drm t}{t}\\
&\lesssim\|g\|_{\CC^\beta}\int_0^1\left(\frac{ts}{(t+s)^2}\right)^{\frac{r}{2}}t^{\frac{\beta}{2}}\|P_tf\|_{L_x^2}\frac{\drm t}{t}.
\end{align*}
This yields
\begin{align*}
\bigg\|s^{-\frac{\alpha+\beta}{2}}\Big\|\int_0^1 Q_sQ_t^{1\bullet}&\left(P_tf\cdot Q_t^2g\right)\frac{\drm t}{t}\Big\|_{L_x^2}\bigg\|_{L^2(s^{-1}\drm s)}\\
&\lesssim\|g\|_{\CC^\beta}\Big\|s^{-\frac{\alpha+\beta}{2}}\int_0^1\left(\frac{ts}{(t+s)^2}\right)^{\frac{r}{2}}t^{\frac{\beta}{2}}\|P_tf\|_{L_x^2}\frac{\drm t}{t}\Big\|_{L^2(s^{-1}\drm s)}\\
&\lesssim\|g\|_{\CC^\beta}\Big\|s^{-\frac{\alpha}{2}}\|P_sf\|_{L_x^2}\Big\|_{L^2(s^{-1}\drm s)}\\
&\lesssim\|f\|_{\CH^\alpha}\|g\|_{\CC^\beta}
\end{align*}
where we used that $\alpha<0$ since $P$ can encode no cancellation and this complete the proof for the third estimate. The proofs for the other estimates on $\P_fg$ are similar and we only give the details for the resonant term. Let $Q\in\mathsf{StGC}^r$ with $r>|\alpha+\beta|$ and recall that $\PI(f,g)$ is a linear combination of terms
$$
\int_0^1P_t^\bullet\left(Q_t^1f\cdot Q_t^2g\right)\frac{\drm t}{t}
$$
with $Q^1,Q^2\in\mathsf{StGC}^{\frac{b}{2}}$ and $P\in\mathsf{StGC}^{[0,b]}$. Given $s\in(0,1]$, we have
$$
\left\|\int_0^1 Q_sP_t^\bullet\left(Q_t^1f\cdot Q_t^2g\right)\frac{\drm t}{t}\right\|_{L_x^2}\lesssim\int_0^s\|Q_t^1f\cdot Q_t^2g\|_{L_x^2}\frac{\drm t}{t}+\int_s^1\left(\frac{s}{t}\right)^{\frac{r}{2}}\left\|Q_t^1f\cdot Q_t^2g\right\|_{L_x^2}\frac{\drm t}{t}
$$
and the result follows again from the lemmas using that $\alpha+\beta>0$.
\end{proof}

\medskip

The dependance of $\PT^s$ with respect to $s$ in given in the following proposition.

\medskip

\begin{proposition}\label{ConvergeEta}
Let $s\in(0,1)$ and a regularity exponent $\beta\in(0,1)$. Given $g\in\CC^\beta$, we have
$$
\|f\mapsto\PT_f^s g\|_{L^2\to\CH^\gamma}\lesssim\frac{s^{\frac{\beta-\gamma}{4}}}{\beta-\gamma}\|g\|_{\CC^\beta}
$$
for any $\gamma\in[0,\beta)$.
\end{proposition}

\medskip

\begin{proof}
Given $f\in L^2$ and $\gamma\in[0,\beta)$, we want to bound the $\CH^\gamma$ norm of $\PT_f^sg$ hence let $Q\in\mathsf{StGC}^r$ with $r>|\gamma|$. Recall that $\PT_f^sg$ is a linear combination of terms of the form
$$
\int_0^s\widetilde Q_t^{1\bullet}\left(P_tf\cdot\widetilde Q_t^2g\right)\frac{\drm t}{t}
$$
with $\widetilde Q^1\in\mathsf{GC}^{\frac{b}{2}-2},\widetilde Q^2\in\mathsf{StGC}^{\frac{b}{2}}$ and $P\in\mathsf{StGC}^{[0,b]}$. Given $u\in(0,1]$, we have
\begin{align*}
\left\|\int_0^sQ_uQ_t^{1\bullet}\left(P_tf\cdot Q_t^2g\right)\frac{\drm t}{t}\right\|_{L_x^2}&\lesssim\int_0^s\left(\frac{tu}{(t+u)^2}\right)^{\frac{r}{2}}\left\|P_tf\cdot Q_t^2g\right\|_{L_x^2}\frac{\drm t}{t}\\
&\lesssim\|g\|_{\CC^\beta}\int_0^s\left(\frac{tu}{(t+u)^2}\right)^{\frac{r}{2}}t^{\frac{\beta}{2}}\|P_tf\|_{L_x^2}\frac{\drm t}{t}.
\end{align*}
This yields
\begin{align*}
\bigg\|u^{-\frac{s}{2}}\Big\|\int_0^sQ_uQ_t^{1\bullet}&\left(P_tf\cdot Q_t^2g\right)\frac{\drm t}{t}\Big\|_{L_x^2}\bigg\|_{L^2(u^{-1}\drm u)}\\
&\lesssim\|g\|_{\CC^\beta}\Big\|u^{-\frac{\gamma}{2}}\int_0^s\left(\frac{tu}{(t+u)^2}\right)^{\frac{r}{2}}t^{\frac{\beta}{2}}\|P_tf\|_{L_x^2}\frac{\drm t}{t}\Big\|_{L^2(u^{-1}\drm u)}\\
&\lesssim\|g\|_{\CC^\beta}\frac{4r}{r^2-\gamma^2}s^{\frac{\beta'-\gamma}{2}}\|u^{-\frac{\beta'-\beta}{2}}\|P_uf\|_{L^2}\|_{L^2(u^{-1}\drm u)}\\
&\lesssim\|g\|_{\CC^\beta}\frac{4r}{r^2-\gamma^2}s^{\frac{\beta'-\gamma}{2}}\frac{2}{k+\beta-\beta'}\|f\|_{\CH^{\beta'-\beta}}\\
&\lesssim\frac{\|g\|_{\CC^\beta}}{1-\beta}\frac{s^{\frac{\beta'-\gamma}{2}}}{k+\beta-\beta'}\|f\|_{\CH^{\beta'-\beta}}
\end{align*}
for any $\beta'\in(\gamma,\beta)$ and $P\in\mathsf{StGC}^k$ using that $r\ge1$. For $k\ge1$, one can take $\beta'=\beta$ and get
$$
\bigg\|u^{-\frac{\gamma}{2}}\Big\|\int_0^sQ_uQ_t^{1\bullet}\left(P_tf\cdot Q_t^2g\right)\frac{\drm t}{t}\Big\|_{L_x^2}\bigg\|_{L^2(u^{-1}\drm u)}\lesssim\frac{s^{\frac{\beta-\gamma}{2}}}{1-\beta}\|g\|_{\CC^\beta}\|f\|_{L^2}.
$$
For $k=0$, we have
\begin{align*}
\bigg\|u^{-\frac{\gamma}{2}}\Big\|\int_0^sQ_uQ_t^{1\bullet}\left(P_tf\cdot Q_t^2g\right)\frac{\drm t}{t}\Big\|_{L_x^2}\bigg\|_{L^2(u^{-1}\drm u)}&\lesssim\frac{\|g\|_{\CC^\beta}}{1-\beta}\frac{s^{\frac{\beta'-\gamma}{2}}}{\beta-\beta'}\|f\|_{L^2}
\end{align*}
hence taking $\beta'=\frac{\gamma+\beta}{2}$ yields
$$
\bigg\|u^{-\frac{\gamma}{2}}\Big\|\int_0^sQ_uQ_t^{1\bullet}\left(P_tf\cdot Q_t^2g\right)\frac{\drm t}{t}\Big\|_{L_x^2}\bigg\|_{L^2(u^{-1}\drm u)}\lesssim\frac{s^{\frac{\beta-\gamma}{4}}}{(1-\beta)(\beta-\gamma)}\|g\|_{\CC^\beta}\|f\|_{L^2}.
$$
\end{proof}

\medskip

\begin{proposition}\label{DivergeEta}
Let $s\in(0,1)$ and a regularity exponent $\beta<2$. Given $g\in\CC^\beta$, we have
$$
\|(\PT_f-\PT_f^s)g\|_{\CH^2}\lesssim s^{\frac{\beta-2}{2}}\|f\|_{L^2}\|g\|_{\CC^\beta}
$$
for any $f\in L^2$.
\end{proposition}

\medskip

\begin{proof}
Given $f\in L^2$, we want to bound the $\CH^2$ norm of $(\PT_f-\PT_f^s)g$ hence let $Q\in\mathsf{StGC}^r$ with $r>2$. It is a linear combination of terms
$$
\int_s^1\widetilde Q_t^{1\bullet}\left(P_tf\cdot\widetilde Q_t^2g\right)\frac{\drm t}{t}
$$
with $\widetilde Q^1\in\mathsf{GC}^{\frac{b}{2}-2},\widetilde Q^2\in\mathsf{StGC}^{\frac{b}{2}}$ and $P\in\mathsf{StGC}^{[0,b]}$. Given $u\in(0,1]$, we have
\begin{align*}
\left\|\int_s^1Q_uQ_t^{1\bullet}\left(P_tf\cdot Q_t^2g\right)\frac{\drm t}{t}\right\|_{L_x^2}&\lesssim\int_s^1\left(\frac{tu}{(t+u)^2}\right)^{\frac{r}{2}}\left\|P_tf\cdot Q_t^2g\right\|_{L_x^2}\frac{\drm t}{t}\\
&\lesssim\|f\|_{L^2}\|g\|_{\CC^\beta}\int_s^1\left(\frac{tu}{(t+u)^2}\right)^{\frac{r}{2}}t^{\frac{\beta}{2}}\frac{\drm t}{t}
\end{align*}
using that $\|P_tf\|_{L^2}\lesssim\|f\|_{L^2}$. This yields
\begin{align*}
\bigg\|u^{-1}\Big\|\int_s^1Q_uQ_t^{1\bullet}&\left(P_tf\cdot Q_t^2g\right)\frac{\drm t}{t}\Big\|_{L_x^2}\bigg\|_{L^2(u^{-1}\drm u)}\\
&\lesssim\|g\|_{\CC^\beta}\Big\|u^{-1}\int_s^1\left(\frac{tu}{(t+u)^2}\right)^{\frac{r}{2}}t^{\frac{\beta}{2}}\frac{\drm t}{t}\Big\|_{L^2(u^{-1}\drm u)}\\
&\lesssim s^{\frac{\beta-2}{2}}\|f\|_{L^2}\|g\|_{\CC^\beta}
\end{align*}
and the proof is complete.
\end{proof}

\medskip

\begin{proposition}\label{CorrectorProof}
Let $\alpha_1\in(0,1)$ and $\alpha_2,\beta\in\IR$. If
$$
\alpha_2+\beta<0\quad\text{and}\quad\alpha_1+\alpha_2+\beta>0,
$$
then $(a_1,a_2,b)\mapsto\DC(a_1,a_2,b)$ extends in a unique bilinear operator from $\CC^{\alpha_1}\times\CC^{\alpha_2}\times\CC^\beta$ to $\CC^{\alpha_1+\alpha_2+\beta}$ and from $\CH^{\alpha_1}\times\CC^{\alpha_2}\times\CC^\beta$ to $\CH^{\alpha_1+\alpha_2+\beta}$.
\end{proposition}

\medskip

\begin{proof}
We first consider $(a_1,a_2,b)\in\CC^{\alpha_1}\times\CC^{\alpha_2}\times\CC^\beta$. We want to compute the regularity of 
$$
\DC(a_1,a_2,b)=\PI\big(\PT_{a_1}a_2,b\big)-a_1\PI(a_2,b)
$$ 
using a family $Q$ of $\mathsf{StGC}^r$ with $r>|\alpha_1+\alpha_2+\beta|$. Recall that a term $\PI(a,b)$ can be written as a linear combination of terms of the form
$$
\int_0^1 P_t^{1\bullet}(Q_t^1a\cdot Q_t^2b)\frac{dt}{t},
$$
while $\PT_ba$ is a linear combination of terms of the form
$$
\int_0^1\widetilde Q_t^{3\bullet}\big(P_t^2b\cdot\widetilde Q_t^4a\big)\frac{dt}{t}
$$
with $Q^1,Q^2,\widetilde Q^4\in\mathsf{StGC}^{\frac{b}{2}}$, $\widetilde Q^3\in\mathsf{GC}^{\frac{b}{2}-2}$ and $P^1,P^2\in\mathsf{StGC}^{[0,b]}$. For $P^2\in\mathsf{StGC}^{[1,b]}$, we already have the correct regularity since
\begin{align*}
\int_0^1\int_0^1Q_uP_t^{1\bullet}&\left(Q_t^1\widetilde Q_s^{3\bullet}\left(P_s^2a_1\cdot\widetilde Q_s^4a_2\right)\cdot Q_t^2b\right)\frac{ds}{s}\frac{dt}{t}  \\
&\lesssim\|a_1\|_{\alpha_1}\|a_2\|_{\alpha_2}\|b\|_\beta\int_0^1\int_0^1\left(\frac{ut}{(t+u)^2}\right)^{\frac{r}{2}}\left(\frac{ts}{(s+t)^2}\right)^{\frac{b}{2}}s^{\frac{\alpha_1+\alpha_2}{2}}t^{\frac{\beta}{2}}\frac{ds}{s}\frac{dt}{t}  \\
&\lesssim\|a_1\|_{\alpha_1}\|a_2\|_{\alpha_2}\|b\|_\beta\ u^{\frac{\alpha_1+\alpha_2+\beta}{2}}
\end{align*}
using that $\alpha_1\in(0,1)$. We consider $P^2\in\mathsf{StGC}^0$ for the remainder of the proof. For all $x\in M$, we have
\begin{align*}
\DC\big(a_1,a_2,b\big)(x)&=\PI\left(\PT_{a_1}a_2,b\right)(x)-a_1(x)\cdot\PI(a_2,b)(x)\\
&=\PI\left(\PT_{a_1}a_2-a_1(x)\cdot a_2,b\right)(x)\\
&\simeq\PI\left(\PT_{a_1-a_1(x)}a_2,b\right)(x),
\end{align*}
since $\PI$ is bilinear and $a_1(x)$ is a scalar and $\PT_1a_1=a_1$ up to smooth terms. Thus we only have to consider a linear combination of terms of the form
$$
\int_0^1\int_0^1 P_t^{1\bullet}\bigg(Q_t^1\widetilde Q_s^{3\bullet}\left(\left(P_s^2a_1-a_1(x)\right)\cdot\widetilde Q_s^4a_2\right)\cdot Q_t^2b\bigg)(x)\,\frac{ds}{s}\frac{dt}{t}
$$
using that $\dst\int_0^1\widetilde Q_s^{3\bullet}\widetilde Q_s^4\frac{ds}{s}=\textup{Id}$ up to smooth terms. This gives $\big(Q_u\DC(a_1,a_2,b)\big)(x)$ as a linear combination of terms of the form
\begin{small}\begin{align*}
&\int K_{Q_u}(x,x')P_t^{1\bullet}\bigg(Q_t^1\widetilde Q_s^{3\bullet}\left(\left(P_s^2a_1-a_1(x')\right)\cdot\widetilde Q_s^4a_2\right)\cdot Q_t^2b\bigg)(x')\,\frac{ds}{s}\frac{dt}{t}\nu(dx')   \\
&= \int K_{Q_u}(x,x')K_{P_t^{1\bullet}}(x',x'')\bigg(Q_t^1\widetilde Q_s^{3\bullet}\left(\left(P_s^2a_1-a_1(x'')\right)\cdot\widetilde Q_s^4a_2\right)\cdot Q_t^2b\bigg)(x'')\frac{ds}{s}\frac{dt}{t}\nu(dx')\nu(dx'')   \\
&\quad +\int\int_0^u K_{Q_u}(x,x')K_{P_t^{1\bullet}}(x',x'')\Big(a_1(x'')-a_1(x')\Big)\left(Q_t^1a_2\cdot Q_t^2b\right)(x'')\frac{dt}{t}\nu(dx')\nu(dx'')   \\
&\quad +\int\int_u^1 K_{Q_u}(x,x')K_{P_t^{1\bullet}}(x',x'')\Big(a_1(x'')-a_1(x')\Big)\left(Q_t^1a_2\cdot Q_t^2b\right)(x'')\frac{dt}{t}\nu(dx')\nu(dx'')   \\
&=: A+B+C.
\end{align*}\end{small}
The term $A$ is bounded using cancellations properties. We have
\begin{align*}
|A|&= \int K_{Q_uP_t^{1\bullet}}(x,x')\bigg(Q_t^1\widetilde Q_s^{3\bullet}\left(\left(P_s^2a_1-a_1(x')\right)\cdot\widetilde Q_s^4a_2\right)\cdot Q_t^2b\bigg)(x') \frac{ds}{s}\frac{dt}{t}\nu(dx')   \\
&\lesssim \|a_1\|_{\alpha_1}\|a_2\|_{\alpha_2}\|b\|_\beta\left(\int_0^u\int_0^1\left(\frac{st}{(s+t)^2}\right)^{\frac{b}{2}}(s+t)^{\frac{\alpha_1}{2}}s^{\frac{\alpha_2}{2}}t^{\frac{\beta}{2}}\frac{ds}{s}\frac{dt}{t}\right.\\
&\qquad \left.+\int_u^1\int_0^1\left(\frac{tu}{(t+u)^2}\right)^{\frac{r}{2}}\left(\frac{st}{(s+t)^2}\right)^{\frac{b}{2}}(s+t)^{\frac{\alpha_1}{2}}s^{\frac{\alpha_2}{2}}t^{\frac{\beta}{2}}\frac{ds}{s}\frac{dt}{t}\right)\\
&\lesssim \|a_1\|_{\alpha_1}\|a_2\|_{\alpha_2}\|b\|_\beta\ u^{\frac{\alpha_1+\alpha_2+\beta}{2}},
\end{align*}
using that $\alpha_1\in(0,1),P^2\in\mathsf{StGC}^0$ and $(\alpha_1+\alpha_2+\beta)>0$. 

\medskip

For the term $B$, we have
\begin{align*}
|B|&\lesssim\|a_1\|_{\alpha_1}\|a_2\|_{\alpha_2}\|b\|_\beta\int_{x',x''}\int_0^uK_{Q_u}(x,x')K_{P_t^{1\bullet}}(x',x'')d(x',x'')^{\alpha_1}t^{\frac{\alpha_2+\beta}{2}}\frac{dt}{t}\nu(dx')\nu(dx'')\\
&\lesssim\|a_1\|_{\alpha_1}\|a_2\|_{\alpha_2}\|b\|_\beta\int_0^ut^{\frac{\alpha_1+\alpha_2+\beta}{2}}\frac{dt}{t}\\
&\lesssim\|a_1\|_{\alpha_1}\|a_2\|_{\alpha_2}\|b\|_\beta\ u^{\frac{\alpha_1+\alpha_2+\beta}{2}},
\end{align*}
using again that $\alpha_1\in(0,1)$ and $(\alpha_1+\alpha_2+\beta)>0$. 

\medskip

Finally for $C$, we also use cancellation properties to get
\begin{align*}
|C|&\lesssim\|a_1\|_{\alpha_1}\|a_2\|_{\alpha_2}\|b\|_\beta \bigg\{\int_{x',x''}\int_u^1K_{Q_u}(x,x')K_{P_t^{1\bullet}}(x',x'')\Big|a_1(x)-a_1(x')\Big|t^{\frac{\alpha_2+\beta}{2}}\frac{dt}{t}\nu(dx')\nu(dx'')   \\
&\quad+ \int_{x',x''}\int_u^1K_{Q_u}(x,x')K_{P_t^{1\bullet}}(x',x'')\Big|a_1(x')-a_1(x'')\Big|t^{\frac{\alpha_2+\beta}{2}}\frac{dt}{t}\nu(dx')\nu(dx'')\bigg\}   \\
&\lesssim\|a_1\|_{\alpha_1}\|a_2\|_{\alpha_2}\|b\|_\beta \bigg\{\int_{x',x''}\int_u^1K_{Q_u}(x,x')K_{P_t^{1\bullet}}(x',x'')d(x,x')^{\alpha_1}t^{\frac{\alpha_2+\beta}{2}}\frac{dt}{t}\nu(dx')\nu(dx'')   \\
&\quad+ \int_{x',x''}\int_u^1K_{Q_u}(x,x')K_{P_t^{1\bullet}}(x',x'')d(x',x'')^{\alpha_1}t^{\frac{\alpha_2+\beta}{2}}\frac{dt}{t}\nu(dx')\nu(dx'')\bigg\}   \\
&\lesssim\|a_1\|_{\alpha_1}\|a_2\|_{\alpha_2}\|b\|_\beta \bigg\{u^{\frac{\alpha_1}{2}}\int_u^1t^{\frac{\alpha_2+\beta}{2}}\frac{dt}{t}+\int_u^1 \left(\frac{tu}{(t+u)^2}\right)^{\frac{r}{2}}t^{\frac{\alpha_1+\alpha_2+\beta}{2}}\frac{dt}{t}\bigg\}  \\
&\lesssim\|a_1\|_{\alpha_1}\|a_2\|_{\alpha_2}\|b\|_\beta\ u^{\frac{\alpha_1+\alpha_2+\beta}{2}},
\end{align*}
using that $\alpha_1\in(0,1)$ and $(\alpha_2+\beta)<0$. In the end, we have
$$
\Big\|Q_u\DC(a_1,a_2,b)\Big\|_\infty\lesssim\|a_1\|_{\alpha_1}\|a_2\|_{\alpha_2}\|b\|_\beta\ u^{\frac{\alpha_1+\alpha_2+\beta}{2}}
$$
uniformly in $u\in(0,1]$, so the proof is complete for $\DC$. The adaptation of the proof to the case $a_1\in\CH^{\alpha_1}$ is left to the reader and follows from the estimates of the Appendix \ref{AppendixApprox}.
\end{proof}

\bigskip

\noindent \textcolor{gray}{$\bullet$} A. Mouzard --  Univ. Rennes, CNRS, IRMAR - UMR 6625, F-35000 Rennes, France   \\
{\it E-mail}: antoine.mouzard@ens-rennes.fr

\end{document}